\documentclass[11pt]{amsart}
\usepackage{amssymb}
\usepackage{amssymb, amsmath}
\usepackage{verbatim}
\usepackage{color}
\let\d=\partial
\let\eps=\varepsilon
\let\wt=\widetilde

\renewcommand{\d}{\partial}

\newcommand{\esp}[1]{\quad \text{#1} \quad}

\renewcommand{\div}{{\rm div\,}}

\newcommand{\cA}{\mathcal{A}}
\newcommand{\cC}{\mathcal{C}}

\newcommand{\cI}{\mathcal{I}}
\newcommand{\cL}{\mathcal{L}}
\newcommand{\cM}{\mathcal{M}}

\newcommand{\cP}{\mathcal{P}}
\newcommand{\cQ}{\mathcal{Q}}
\newcommand{\cR}{\mathcal{R}}

\newcommand{\ca}{\check{a}}
\newcommand{\cu}{\check{u}}
\newcommand{\cw}{\check{w}}
\newcommand{\crho}{\check{\rho}}

\def\da{\delta\!a}

\newcommand{\N}{\mathbb N}

\newcommand{\R}{\mathbb R}
\newcommand{\Z}{\mathbb Z}

\newcommand{\ddj}{\Dot{\Delta}_j}

\newcommand{\with}{\quad\!\hbox{with}\!\quad}
\newcommand{\andf}{\quad\!\hbox{and}\!\quad}

\newtheorem{thm}{Theorem}[section]
\newtheorem{rmq}{Remark}[section]

\title[High viscosity limit]{
High viscosity limit for the multi-dimensional compressible
Navier-Stokes equations}
\author{Rapha\"el Danchin\\ Universit\'e Paris-Est Cr\'eteil}

\address{Univ Paris Est Creteil, Univ Gustave Eiffel, CNRS, LAMA UMR8050, F-94010 Creteil, France \&
 Institut Universitaire de France}
 \email{raphael.danchin@u-pec.fr}

\date\today

\begin{document}
    
    \begin{abstract} We investigate the high viscosity limit (also called inertial limit) 
    of the barotropic compressible Navier-Stokes equations
    supplemented with initial data which are perturbations of a stable constant solution. 
    In the case of constant viscosity coefficients, we establish that, after diffusive rescaling, the density  tends to 
    satisfy a transport equation with nonlinear damping  which is globally well-posed, even for large data.
        Similar results are proved  for variable viscosity coefficients. In this latter case, the damping term in the
         limit equation of the density is nonlocal.
         \end{abstract}
    \maketitle
    We consider the barotropic compressible Navier-Stokes equations
     $$   \left\{\begin{aligned}
    &\d_t \rho+\div( \rho u)=0,\\
   &\rho(\d_tu+u\cdot\nabla u) +\nu  \cA u+c^2\nabla(P(\rho))=0\end{aligned}\right. \leqno(CNS_{\nu,c})$$
   where  $u=u(t,x)\in\R^d$ stands for the velocity field of the compressible fluid    under consideration and   
   $\rho=\rho(t,x)\in\R_+,$ for its (renormalized) density, tending to $1$ at infinity.
The function $P$ is the (renormalized) pressure, supposed smooth and such that $P(1)=0$ and $P'(1)=1.$
The positive real number  $c$ stands for the sound speed  at reference  density $\rho=1.$
   Finally, $\nu>0$ is the total kinematic viscosity, and  we put
   \begin{equation}\label{def:A}   \cA :=-\mu\Delta -(1-\mu)\nabla\div= -\Delta \cQ  - \mu \Delta \cP, \end{equation}
    where $\cP$ and $\cQ$ are the  Helmholtz projectors on solenoidal and potential vector-fields, and $\mu>0$ is the relative shear viscosity.
      \smallbreak
   In the last section of our recent work \cite{D25} 
   dedicated to the one-dimensional case in mass Lagrangian coordinates, we pointed out
   an intriguing asymptotics for $\nu$ going to $\infty$: the density tends to the solution of 
   a nonlinear Ordinary Differential Equation which  has global solutions whenever the pressure law is nondecreasing.
   In the present work, we investigate whether similar results are valid in the multi-dimensional case. 
   \medbreak
   Before going into more details, let us give  a few benchmarks in the study of the initial value problem for 
   the compressible Navier-Stokes equations  in the multi-dimensional case.
For smooth initial data, the proof of existence of classical solutions on a small time interval 
   goes back to the seminal works by Serrin \cite{Serrin} and Nash \cite{Nash} in the early sixties. 
   Twenty years later,  Matsumura and Nishida in \cite{MN} established the global well-posedness of
   the compressible Navier-Stokes equations 
   in the case where the initial data $(\rho_0,u_0)$ are Sobolev perturbation (with sufficiently high regularity index) 
   of some constant state $(\rho,u)=(\bar\rho,0)$ such that $\bar\rho>0$ and $P'(\bar\rho)>0.$ 
   In the 90ies P.-L. Lions in \cite{PLL} constructed global solutions \emph{\`a la Leray} 
   satisfying the energy inequality associated to $(CNS_{\nu,c})$  for general finite energy data 
  in the so-called isentropic case $P(\rho)=a\rho^\gamma$ with $a>0$ and $\gamma>1$ large enough
   (see also    the work by Feireisl  \cite{Feireisl}),
   while Hoff in \cite{Hoff95} obtained global solutions with `intermediate' regularity.

   In all the aforementioned  pioneering works however,  due to the use of a nonhomogeneous functional framework,    
   tracking the dependence of the solutions with respect to the total viscosity is delicate, and 
    it is hard to justify the asymptotics $\nu$ going to $\infty$. 
    In fact,  the only works we are aware of is the very recent paper by C. Yu \cite{Yu} where 
    the weak convergence (up to subsequence)  to the same limit system as the one 
    which will be presented below is justified, for the finite energy solutions, 
    and our recent work \cite{D25} already mentioned above, which  is specific to the one-dimensional 
    situation.

   The present work aims at filling the gap in the multi-dimensional case, for small data with critical regularity.
   We want to obtain results of strong convergence, with algebraic decay rates with respect to the 
   total viscosity $\nu.$

   \section{Results}
   
   An essential ingredient for studying the high viscosity limit  is to have uniform estimates. This motivates us to use a homogeneous functional framework as it allows estimates in the general case to be deduced from the case $\nu=1$.
  We here adopt the homogeneous (or rather quasi-homogeneous) 
   critical functional  framework that has been proposed by the author in \cite{D00}. 
  \smallbreak
  To be more specific, introducing some notation is in order. 
Let $(\ddj)_{j\in\Z}$ be the family of spectral truncation operators defined in \cite[Chap. 2]{BCD}.
One can set for instance $\ddj:=\varphi(2^{-j}D)$ with $\varphi(\xi):=\chi(\xi/2)-\chi(\xi)$ 
and  $\chi\in\cC_c^\infty(B(0,4/3))$ nonnegative and with value $1$
  on $B(0,3/4).$

  For all $s\in\R,$ we define the  homogeneous Besov space $\dot B^s_{2,1}$ 
  to be the set of tempered distributions $u$ going weakly to $0$
  at infinity  (that is, $\lim \|\chi(2^{-j}D)u\|_{L^\infty}$ tends to $0$ when $j$ goes to $-\infty$), 
  and such that 
   $$  \|u\|_{\dot B^s_{2,1}}:=\sum_{j\in\Z} 2^{js}\|\ddj u\|_{L^2}<\infty.$$
  We need also spaces with different regularity indices for low and high frequencies. 
  To do this,   we set for $\alpha\in\R_+$ and $\sigma\in\R,$
\begin{eqnarray}\label{eq:not1}
&\displaystyle\|z\|_{\dot B^\sigma_{2,1}}^{\ell,\alpha}:= \sum_{j\leq 1+\log_2 \alpha} 2^{j\sigma}\|\ddj z\|_{L^2}\andf
\|z\|_{\dot B^\sigma_{2,1}}^{h,\alpha}:= \sum_{j\geq \log_2 \alpha} 2^{j\sigma}\|\ddj z\|_{L^2},\\
\label{eq:not2}
&\displaystyle z^{\ell,\alpha}:= \sum_{j\leq \log_2 \alpha} \ddj z\andf
z^{h,\alpha}:= \sum_{j>\log_2 \alpha} \ddj z.
\end{eqnarray}
The intentional redundancy between low and high frequencies in~\eqref{eq:not1} will  enable us to use that
\begin{equation}\label{eq:redond}
\|z^{\ell,\alpha}\|_{\dot B^s_{2,1}} \leq \|z\|^{\ell,\alpha}_{\dot B^s_{2,1}} \andf
\|z^{h,\alpha}\|_{\dot B^s_{2,1}} \leq \|z\|^{h,\alpha}_{\dot B^s_{2,1}}.
\end{equation}
Performing the change of unknowns 
\begin{equation}\label{eq:rescaling}
\wt a(t,x):={\rho}\biggl(\frac{\nu}{c^2}\,t,\frac{\nu}{c}\,x\biggr)-1
\andf \wt u(t,x):=\frac{u}{c}\biggl(\frac{\nu}{c^2}\,t,\frac{\nu}{c}\,x\biggr),
\end{equation}
transforms the original equations into 
       \begin{equation}  \left\{\begin{aligned}
    &\d_t \wt a+\div( (1+\wt a)\wt u)=0,\\
   &(1+\wt a)(\d_t\wt u+\wt u\cdot\nabla \wt u) +\cA \wt u+P'(1+\wt a)\nabla \wt a=0.\end{aligned}\right. \label{eq:CNS}
   \end{equation}
    The classical  global well-posedness theorem in critical spaces \cite{D00} asserts that there exist two positive numbers
         $\eta_0=\eta_0(d,\mu,P')$ and $C=C(d,\mu,P')$  such that  if 
   \begin{equation}\label{eq:smalldata0}
   \|(\wt a_0,\wt u_0)\|_{\dot B^{\frac d2-1}_{2,1}} + \|\wt a_0\|_{\dot B^{\frac d2}_{2,1}}\leq \eta_0\end{equation}
   then \eqref{eq:CNS} has a unique global solution $(\wt a,\wt u)$ 
   with $\wt a\in\cC_b(\R_+;\dot B^{\frac d2-1}_{2,1}\cap \dot B^{\frac d2}_{2,1})\cap L^2(\R_+;\dot B^{\frac d2}_{2,1})$
   and $\wt u \in\cC_b(\R_+;\dot B^{\frac d2-1}_{2,1})\cap L^1(\R_+;\dot B^{\frac d2+1}_{2,1})$
      verifying for all $t\in\R_+,$
   \begin{multline}\label{eq:globcns}  \|(\wt a,\wt u)\|_{L^\infty_t(\dot B^{\frac d2-1}_{2,1})} + \|\wt a\|_{L^\infty_t(\dot B^{\frac d2}_{2,1})}
   +\|\nabla \wt u\|_{L^1_t(\dot B^{\frac d2}_{2,1})}+ \|(\wt a,\wt u)\|_{L^2_t(\dot B^{\frac d2}_{2,1})}\\
   + \|\nabla \wt a\|^{\ell,1}_{L^1_t(\dot B^{\frac d2}_{2,1})}+ \|\wt a\|^{h,1}_{L^1_t(\dot B^{\frac d2}_{2,1})}
   \leq C\bigl( \|(\wt a_0,\wt u_0)\|_{\dot B^{\frac d2-1}_{2,1}} + \|\wt a_0\|_{\dot B^{\frac d2}_{2,1}}\bigr)\cdotp
   \end{multline}
   Recall  (see e.g. \cite[Chap. 2]{BCD}) that  for all regularity index $s,$ frequency threshold $\mathfrak f$ and positive real number $\lambda,$ we have 
   \begin{equation}\label{eq:rescalingb}
\|z(\lambda\cdot)\|_{\dot B^s_{2,1}} \simeq\lambda^{s-\frac d2}\|z\|_{\dot B^s_{2,1}},\quad\!\! 
\|z(\lambda\cdot)\|^{\ell,\mathfrak f}_{\dot B^s_{2,1}} \simeq\lambda^{s-\frac d2}\|z\|_{\dot B^s_{2,1}}^{\ell,\mathfrak f\lambda^{-1}},\quad\!\!
\|z(\lambda\cdot)\|^{h,\mathfrak f}_{\dot B^s_{2,1}} \simeq\lambda^{s-\frac d2}\|z\|_{\dot B^s_{2,1}}^{h,\mathfrak f\lambda^{-1}}.
\end{equation}
Hence, from \eqref{eq:rescaling}, 
 \eqref{eq:smalldata0} and \eqref{eq:rescalingb},  we deduce that, if    
      $$ \nu^{-1}\bigl(c\|a_0\|_{\dot B^{\frac d2-1}_{2,1}} + \|u_0\|_{\dot B^{\frac d2-1}_{2,1}}\bigr) + \|a_0\|_{\dot B^{\frac d2}_{2,1}}\leq \eta_0 $$
   then $(CNS_{\nu,c})$ has a unique global solution $(\rho=1+a,u)$ such that
       \begin{multline}\label{eq:uniform1}  \nu^{-1}\bigl(c\|a\|_{L^\infty_t(\dot B^{\frac d2-1}_{2,1})} 
       + \|u\|_{L^\infty_t(\dot B^{\frac d2-1}_{2,1})}\bigr) 
      + \|a\|_{L^\infty_t(\dot B^{\frac d2}_{2,1})}      + c^{3/2}\nu^{-1/2}\|a\|_{L^2_t(\dot B^{\frac d2}_{2,1})}\\+ \nu^{-1/2}\|u\|_{L^2_t(\dot B^{\frac d2}_{2,1})}          +\|\nabla u\|_{L^1_t(\dot B^{\frac d2}_{2,1})}    
   + {c}\|\nabla a\|^{\ell,c\nu^{-1}}_{L^1_t(\dot B^{\frac d2}_{2,1})}\\+ c^2\nu^{-1}\|a\|^{h,c\nu^{-1}}_{L^1_t(\dot B^{\frac d2}_{2,1})}
   \leq C\bigl(\nu^{-1}(c\|a_0\|_{\dot B^{\frac d2-1}_{2,1}} + \|u_0\|_{\dot B^{\frac d2-1}_{2,1}}) 
     + \|a_0\|_{\dot B^{\frac d2}_{2,1}}\bigr)\cdotp
   \end{multline}
 
Inequality \eqref{eq:uniform1}  ensures that for fixed data,  $a$  and $\nabla u$ are uniformly bounded in $L^\infty(\R_+;\dot B^{\frac d2}_{2,1})$
and $L^1(\R_+;\dot B^{\frac d2}_{2,1}),$ respectively. Hence, up to subsequence,
$a$  (resp. $\nabla u$) converges  weakly * 
 to some $b$ in $L^\infty(\mathbb R_+;\dot B^{\frac d2}_{2,1})$ (resp. $\nabla v$ in $\cM(\mathbb R_+;\dot B^{\frac d2}_{2,1})$).
 However, without stronger assumptions, it  looks difficult to have more accurate information. 
In fact, a rough analysis of the linearized equations (see e.g. the appendix of \cite{D25}) reveals that for $\nu$ going to infinity, we have
$a(t)\simeq e^{-\frac t\nu}a_0$ which means that, contrary to what one might think, having a large viscosity 
\emph{decreases} the convergence speed  to some constant state when $t$ goes to infinity. 
This is sometimes referred to as the phenomenon of overdamping, leading to a sluggish convergence when the viscosity increases. 
\smallbreak
In order to get more insight, it is wise to look at time scale of order $\nu,$ performing the following \emph{diffusive} change of variables:
 \begin{equation}\label{eq:diffusive}(\rho,u)(t,x)= (\check \rho, \nu^{-1}\check u)(\nu^{-1}t, x)\esp{and}    (\rho_0,u_0)(x)= (\check \rho_0,\nu^{-1}\check u_0)(x).
 \end{equation} 
 Setting   \begin{equation}\label{eq:Q}\check\rho=1+\ca\andf  
 Q(z):= P(1+z),\end{equation} we see that $(\rho,u)$ satisfies  
  $(CNS_{\nu,c})$  if and only if $(\check a,\check u)$ is a solution of
  \begin{equation}\label{eq:cnsdiff}   \left\{\begin{aligned}
    &\d_t \check a+\div( (1+\check a)\check u)=0,\\   
    &\nu^{-2}(1+\check a)\bigl(\d_t\check u+\check u\cdot\nabla \check u\bigr) +\cA \check u+c^2\nabla(Q(\check a))=0.\end{aligned}\right.\end{equation}
       Therefore, since $\cA(-\Delta)^{-1}\nabla=\nabla,$
        one can expect the following `effective velocity':  
    \begin{equation}\label{def:w}
     \check w:=\check u+c^2(-\Delta)^{-1}\nabla (Q(\check a))\end{equation}  
     to   tend to $0$ for $\nu$ going to $\infty.$
    Now, replacing $\check u$ by $\check w$ in the equation for $\check a,$ we get 
    \begin{equation}\label{eq:rescaledmass}\partial_t\check a
    -c^2\div\bigl((1\!+\!\check a)(-\Delta)^{-1}\nabla (Q(\check a))\bigr)=-\div((1\!+\!\check a )\check w).\end{equation}
   Hence, if it is true that $\check w$  tends to $0$ in some reasonable space, then $\ca$ tends to the solution $b$ of the following damped 
   transport equation:
   \begin{equation}\label{eq:limit}\partial_tb +v\cdot\nabla b+ c^2(1+b)Q(b) =0\with v:=- c^2(-\Delta)^{-1}\nabla(Q(b)).\end{equation}
   The main goal of the present  paper is to justify this heuristics rigorously. 
   Our results are summarized in the following statement.
   \begin{thm}\label{thm:main}
   There exists a positive real number $\eta_0=\eta_0(d,\mu,P')$ such that for any family
    $(\rho_0^\nu,u_0^\nu)_{\nu>0}$ of data with $a_0^\nu:=\rho_0^\nu-1$ in $\dot B^{\frac d2-1}_{2,1}$ and 
   $u_0^\nu$ in $\dot B^{\frac d2-1}_{2,1}$ satisfying
   \begin{equation}\label{eq:smalldata}
   \nu^{-1}c\|a_0^\nu\|_{\dot B^{\frac d2-1}_{2,1}} +    \|a_0^\nu\|_{\dot B^{\frac d2}_{2,1}}+  \nu^{-1} \|u_0^\nu\|_{\dot B^{\frac d2-1}_{2,1}}
   \leq\eta_0 \end{equation}
  Equations $(CNS_{\nu,c})$ supplemented with initial data $(\rho_0^\nu,u_0^\nu)$ admit a unique
   global solution $(\rho^\nu=1+a^\nu,u^\nu)$ such that 
   $$   a^\nu\in \cC_b(\R_+;\dot  B^{\frac d2-1}_{2,1}\cap\dot B^{\frac d2}_{2,1})
   \cap L^2(\R_+;\dot B^{\frac d2}_{2,1})\andf 
    u^\nu\in \cC_b(\R_+;\dot  B^{\frac d2-1}_{2,1})   \cap L^1(\R_+;\dot B^{\frac d2+1}_{2,1}).$$
    Furthermore, Inequality \eqref{eq:uniform1} holds true for some $C=C(d,\mu,P').$ 
    \medbreak
    If, in addition, $(a_0^\nu,u_0^\nu)\in\dot B^s_{2,1}$ for some $s\in(-d/2,d/2-1],$ then 
    $$   a^\nu\in \cC_b(\R_+;\dot  B^{s}_{2,1})\andf 
    u^\nu\in \cC_b(\R_+;\dot  B^{s}_{2,1})   \cap L^2(\R_+;\dot B^{s+1}_{2,1}),$$ 
   and we have, setting $\check\nu:=\nu c^{-1},$
   \begin{multline}\label{eq:uniform0}   c\|a^\nu\|_{L^\infty_t(\dot B^s_{2,1})} +\|u^\nu\|_{L^\infty_t(\dot B^s_{2,1})}
   +\nu^{1/2}\|u^\nu\|_{L^2_t(\dot B^{s+1}_{2,1})}
  +\nu \|u^\nu\|_{L^1_t(\dot B^{s+2}_{2,1})}^{\ell,\check\nu^{-1}}
    +\nu \|w^\nu\|_{L^1_t(\dot B^{s+2}_{2,1})}^{h,\check\nu^{-1}}\\
    +c\nu \|a^\nu\|_{L^1_t(\dot B^{s+2}_{2,1})}^{\ell,\check\nu^{-1}}
   +c^3\nu^{-1} \|a^\nu\|_{L^1_t(\dot B^{s}_{2,1})}^{h,\check\nu^{-1}}
   \leq CY^s_{\nu,c,0} \with {Y^s_{\nu,c,0}:=  c\|a^\nu_0\|_{\dot B^s_{2,1}} +\|u^\nu_0\|_{\dot B^s_{2,1}}},\end{multline}
     where $w^\nu:=u^\nu+c^2\nu^{-1}(-\Delta)^{-1}\nabla(P(\rho^\nu))$ 
     and $C=C(s,d,\mu,P').$
     \medbreak
     Finally, if $Y^s_{\nu,c,0}$ is bounded for $\nu\to\infty$ 
         and  $a_0^\nu\to b_0$ in $\dot B^s_{2,1}$ for some $b_0\in\dot B^s_{2,1}\cap\dot B^{\frac d2}_{2,1}$ then: 
     \begin{enumerate}
     \item Equation \eqref{eq:limit} supplemented with initial data $b_0$ admits a unique global solution
      $b$ in $\cC_b(\R_+;\dot B^{\frac d2}_{2,1}\cap\dot B^s_{2,1}))
     \cap L^1(\R_+;\dot B^{\frac d2}_{2,1}\cap\dot B^s_{2,1})$ 
   satisfying 
   $$   \|b\|_{L^\infty_t(\dot B^\sigma_{2,1})} +\frac{c^2}2    \|b\|_{L^1_t(\dot B^\sigma_{2,1})} 
 \leq \|b_0\|_{\dot B^\sigma_{2,1}}\quad\hbox{for all }\ \sigma\in[s,d/2].$$    
     \item Let  $v:=-c^2(-\Delta)^{-1}\nabla(Q(b))$ and assume that  $s\in[d/2-2,d/2-1].$  
     Denoting $\alpha_s:=\frac d2-1-s,$  $p_s:=(\frac d4-\frac s2)^{-1}$ and $p'_s$ the conjugate exponent,  there holds for all $t\in\R_+,$
  {\begin{align}\|a^\nu(t)-b(\nu^{-1}t)\|_{\dot B^{\frac d2-1}_{2,1}} &\leq C\Bigl(\|a_0^\nu-b_0\|_{\dot B^{\frac d2-1}_{2,1}}
     +\check\nu^{-\alpha_s} \bigl(c^{-1}Y^s_{\nu,c,0}+\|b_0\|_{\dot B^s_{2,1}}\bigr)\Bigr),\label{eq:asympa}\\ 
       \label{eq:asympv}
      \|u^\nu-\nu^{-1}v(\nu^{-1}\cdot)\|_{L_t^{p_s}(\dot B^{\frac d2}_{2,1})} &\leq
       C\Bigl(\nu^{-\frac1{p'_s}}c^{\frac2{q'_s}}\|a_0^\nu-b_0\|_{\dot B^{\frac d2-1}_{2,1}}
       +\nu^{-\frac1{p_s}}\bigl(Y_{\nu,c,0}^s+c\|b_0\|_{\dot B^s_{2,1}}\bigr)\Bigr)\cdotp
    \end{align}}
   \smallbreak  \item If $d\geq3$ and $s=d/2-2$ then $a^\nu-b(\nu^{-1}\cdot)$ converges uniformly to zero 
     on $\R_+\times\R^d.$
     \smallbreak\item If $d=2$ and in addition $(a_0^\nu,b_0^\nu)$ is bounded in the space $\dot B^{-1}_{2,\infty}$ 
     (defined  in \eqref{eq:besovinfini})
     then the regularity $\dot B^{-1}_{2,\infty}$ is propagated for all time, 
     and the uniform convergence  holds true. 
          \end{enumerate}
          \end{thm}
   \begin{rmq}Observe  that $\|\nu^{-1}v(\nu^{-1}\cdot)\|_{L^{p_s}(\R_+;\dot B^{\frac d2}_{2,1})}$ is of order $\nu^{-\frac1{p'_s}}$
   and that $p_s<p'_s$ if $s<-1+d/2.$ 
   Hence, if  $a_0^\nu$ tends to $b_0$ in $\dot B^{\frac d2-1}_{2,1}$ when $\nu$ goes to $\infty,$    then Inequality \eqref{eq:asympv} ensures  that 
   $\nu^{-1}v(\nu^{-1}\cdot)$ is indeed the leading order term of the expansion of $u^\nu$ for $\nu$ going to $\infty.$
   \end{rmq}
    \begin{rmq} The limit equation for the density is globally well-posed for large data. Studying whether
    we can still have   a global solution $(\rho^\nu,u^\nu)$ for large enough $\nu$ without assuming that $a_0^\nu$ 
    is small will be the subject of  a forthcoming work.
   \end{rmq}
   \begin{rmq} The same result holds true for the variable viscosity case, except that 
  we do not know whether we have uniform convergence, nor whether the limit equation is globally well-posed 
   for large data (see more details in Section \ref{s:variable}). 
   \end{rmq}
   \begin{rmq}  The asymptotic behavior  of the solutions if only the bulk  viscosity tends to  infinity
   (that is,  the viscous stress tensor reads $\nu\Delta\cQ+\mu\Delta\cP$ with fixed $\mu$ and $\nu\to\infty$)  is completely different.
    Under appropriate assumptions, the solutions then tend toward those of the incompressible Navier-Stokes equations, which may be non-homogeneous (see \cite{DM-adv,DM-ripped} for more details).  \end{rmq}
  \begin{rmq}
   The importance of the effective velocity defined in \eqref{def:w} and of its divergence 
    (the so-called viscous effective flux)
   in the study of the compressible Navier-Stokes equations has been pointed out 
   much earlier  by  Hoff \cite{Hoff87}, Serre \cite{Serre} and  P.-L. Lions \cite{PLL} 
   as it provides  key regularity and compactness  to construct global weak solutions. 
   It has also been used by  Haspot \cite{Haspot} for improving the global existence result of \cite{D00}. 
\end{rmq}
\subsection*{Notation} Throughout the paper, $C$  denotes harmless positive constants that may change from line to line, and 
$A\lesssim B$ means $A\leq CB.$  We sometimes denote $C(a,b,\dots)$ to emphasize the dependency 
of $C$ with respect to parameters $a,b,\dots$ Furthermore, the letter $k$ or its variations $k_1,$ $k_2$ and so on  designate various smooth functions 
that vanish at $0.$  

Finally, for any Banach space $X$ and Lebesgue exponent $r\in[1,\infty],$
the notation $\|f\|_{L^r_t(X)}$ is used to designate the norm of measurable functions $f$ from $(0,t)$ to $X$ 
such that $t\mapsto \|f(t)\|_X$ is in the Lebesgue space $L^r(0,t).$

   
   \section{The proof of the main result}\label{s:constant}
   
   Let us shortly enumerate the main steps of the proof. As a first, in order to get the desired result of convergence for the density, we perform the rescaling \eqref{eq:diffusive}.
   Denoting by  $(\crho^\nu,\cu^\nu)$ the rescaled solution, then $\ca^\nu:=\crho^\nu-1,$ we introduce
    the effective velocity
    \begin{equation}\label{def:w2}\cw^\nu:=\cu^\nu+c^2(-\Delta)^{-1}\nabla(Q(\ca^\nu)).
    \end{equation}
    The key is to prove that    $\check w^\nu$  tends strongly to zero 
   when $\nu$ goes to infinity, for some suitable norm. 
   To achieve this, the starting point is the counterpart of Inequality \eqref{eq:uniform1} for the scaling under consideration. 
   Although it does not readily ensure that $\check w^\nu$ tends to zero,  it provides us with uniform bounds which will be crucial in the rest of the proof.  
   Step 2 consists in proving a global-in-time control for the  $\dot B^s_{2,1}$ norm of the solution with $s<d/2-1.$   In 
   contrast with the first step,  all the components of the solution are considered in the same space. 
   The main outcome of this step is that, after rescaling, the  part of $\check w^\nu$ corresponding to frequencies higher than $\check\nu^{-1}$
    tends to $0$  in $L^1(\R_+;\dot B^{s+2}_{2,1})$ with the rate $\nu^{-1}.$
   The next step is to show that the limit equation \eqref{eq:limit} is globally well-posed in $\dot B^{\frac d2}_{2,1},$ and  that
   any additional regularity $\dot B^s_{2,1}$ with 
    \begin{equation}\label{eq:s}
   -d/2<s\leq -1+d/2\end{equation}  (and even regularity $\dot B^{-\frac d2}_{2,\infty}$) is propagated. 
   The last two steps are devoted to proving the convergence of $(\check a^\nu,\check u^\nu)$ 
   to $(b,v)$ satisfying \eqref{eq:limit}. 
   As a first, we show convergence for weaker norms (corresponding to a loss of one derivative).
   Then, we point out that this loss does not occur if we consider    
   the Lagrangian coordinates corresponding 
   to the flows of the vector-fields $\check u^\nu$ and $v,$ respectively. 
   As a by-product, we eventually obtain  the    uniform convergence of $\check a^\nu$ to~$b,$ with an explicit rate.

   \subsection{Uniform bounds in critical norm}
  Equations \eqref{eq:cnsdiff} correspond to  $(CNS_{\nu^2,c\nu}).$
  As for $(CNS_{\nu,c})$ , the threshold between low and high frequencies is at $\check\nu^{-1}$ with 
  $\check \nu=\nu c^{-1}$ and,  due to  \eqref{eq:uniform1}, we have
     \begin{multline}\label{eq:uniform2}    \check\nu^{-1}\|\check a^\nu\|_{L^\infty_t(\dot B^{\frac d2-1}_{2,1})} 
       +  \nu^{-2}\|\check u^\nu\|_{L^\infty_t(\dot B^{\frac d2-1}_{2,1})} 
      + \|\check a^\nu\|_{L^\infty_t(\dot B^{\frac d2}_{2,1})}
          +\|\nabla \check u^\nu\|_{L^1_t(\dot B^{\frac d2}_{2,1})}    
     +\|\nabla \check w^\nu\|_{L^1_t(\dot B^{\frac d2}_{2,1})}^{h,\check\nu^{-1}}\\
   + c^2\bigl( \check\nu\|\nabla \check a^\nu\|^{\ell,\check\nu^{-1}}_{L^1_t(\dot B^{\frac d2}_{2,1})}+ \|\check a^\nu\|^{h,\check\nu^{-1}}_{L^1_t(\dot B^{\frac d2}_{2,1})}\bigr)
   \lesssim   \check\nu^{-1}\|\check a_0^\nu\|_{\dot B^{\frac d2-1}_{2,1}} +  \nu^{-2}\|\check u_0^\nu\|_{\dot B^{\frac d2-1}_{2,1}} 
   +\|\check a_0^\nu\|_{\dot B^{\frac d2}_{2,1}}.
   \end{multline}

 \subsection{Control of lower order norms}
  
Granted with   \eqref{eq:uniform2}, we want to glean algebraic decay with respect to $\nu$ for some part of the solution 
(for the high (resp. low) frequencies of $\cw^\nu$ (resp. $\cu^\nu$) in particular). It will be achieved by studying the evolution of regularity $\dot B^s_{2,1}$ of the overall solution, for some $s$ satisfying \eqref{eq:s}.  
The important point is that  $\ca^\nu$ and $\cu^\nu$ are considered at the \emph{same} level of 
 regularity\footnote{For technical reasons that appear in dimension $2,$  we shall also consider the limit
 regularity   $\dot B^{-\frac d2}_{2,\infty}$.}. 
 
 As for proving \eqref{eq:uniform2}, the overall strategy is to  decompose the solution into frequencies lower and higher than $\check\nu^{-1}.$  
 We estimate  the low frequencies of the solution  leveraging the parabolic properties of
the linearized compressible Navier-Stokes equations.
For the high frequencies, we rather consider the pair
  $(\check a^\nu,\check w^\nu)$ (with $\check w^\nu$ defined in \eqref{def:w2}) which satisfies  
  a coupling between a damped transport equation
  for $\check a^\nu$, and a parabolic equation for $\check w^\nu.$
\smallbreak  
 In light of \eqref{eq:rescalingb}  and of  the  rescaling 
\begin{equation}\label{eq:rescaling2}
(\check a^\nu,\check u^\nu)(t,x)=(a, c\nu u)(c^2t,\check \nu^{-1}x),
\end{equation}
it suffices to consider the case $\nu=c=1.$ 
Then, it is known from \cite{D00} that if $\eta_0$ in \eqref{eq:smalldata} is small enough then \eqref{eq:CNS} admits
a unique global solution satisfying \eqref{eq:globcns}.  In what follows this solution will be just denoted
by $(a,u).$
Note that  $a$  remains close to $0$ in the following sense: 
\begin{equation}\label{eq:smalla}
\forall t\in\R_+,\; \|a(t)\|_{\dot B^{\frac d2}_{2,1}}\leq C\eta_0\ll1,\end{equation}
hence also  $\|a(t)\|_{L^\infty}$ is uniformly small with respect to $\nu$ and $t,$ due to the critical embedding
 \begin{equation}\label{emb:crit}
 \dot B^{\frac d2}_{2,1}\hookrightarrow \cC_b.\end{equation}
We skip the rigorous justification of the fact that regularity $\dot B^s_{2,1}$ is indeed propagated, since  
it may be achieved by means of  the construction procedure used in e.g. \cite[Chap. 10]{BCD}. We then concentrate on the proof
of a priori estimates.
\medbreak
As $\nu=c=1,$ we have  to remember that  the modified velocity is given by 
\begin{equation}\label{eq:parabolic}
w=u+(-\Delta)^{-1}\nabla(Q(a)),\end{equation}
and thus 
\begin{equation}\label{eq:cAw}\cA w=\cA u+\nabla (Q(a)).\end{equation}

      \subsubsection*{Estimates of the low frequencies}

In terms of $a,$ $u$ and $w,$  \eqref{eq:CNS} may be rewritten
   \begin{equation}\label{eq:systemlf} \left\{\begin{aligned}
    &\d_t a+\div( a u)+\div u=0,\\
   &\d_t u+u\cdot\nabla  u +\cA u+\nabla a= g:= k_1(a)\,\cA w + \nabla ( ak_2(a))\end{aligned}\right. \end{equation}
where, as said before, $k_1$ and $k_2$  designate generic smooth functions that vanish at $0.$
\medbreak
For all $j\in\Z,$ set $a_j:=\ddj a,$ $u_j:=\ddj u$ and $g_j:=\ddj g.$ Then, applying $\ddj$ to \eqref{eq:systemlf} yields
  $$   \left\{\begin{aligned}
    &\d_t a_j+\div( a_j u)+\div u_j=R_j^1:=\div\bigl([u,\ddj]a\bigr),\\
   &\d_t u_j+u\cdot\nabla  u_j +\cA u_j+\nabla a_j-g_j=  R_j^2:=[u,\ddj]\cdot\nabla u.\end{aligned}\right. $$
   Let $\kappa>0$ be suitably small, and set
   $$
   \cL_j^2:=\|(a_j,u_j)\|_{L^2}^2+2\kappa\int u_j\cdot\nabla a_j.$$
   Direct computations give
   $$\begin{aligned}\frac12\frac d{dt}\|(a_j,u_j)\|_{L^2}^2\!+\!\|\nabla\cQ u_j\|_{L^2}^2\!+\!\mu\|\nabla\cP u_j\|_{L^2}^2
   &=\frac12\int(|u_j|^2-a_j^2)\div u+\!\int\! \bigl(R_j^1 a_j+(R_j^2+g_j)\cdot u_j\bigr)\\
   \andf\!\!    \frac d{dt}\int u_j\cdot\nabla a_j-\|\div u_j\|_{L^2}^2\!+\!\|\nabla a_j\|_{L^2}^2&=-\int\cA u_j\cdot\nabla a_j\\
   +\int \bigl(\div(a_j u)\,\div u_j-(u\cdot&\nabla u_j)\cdot\nabla a_j\bigr) +\int \bigl(u_j\cdot\nabla R_j^1
   +(R_j^2+g_j)\cdot\nabla a_j \bigr)\cdotp\end{aligned}$$
   Performing suitable integrations by parts reveals that 
   $$  \int \bigl(\div(a_j u)\div u_j-(u\cdot\nabla u_j)\cdot\nabla a_j\bigr)=\int \div u \,\div(a_ju_j)
   -\int(u_j\cdot\nabla u)\cdot\nabla a_j.$$
   Hence, observing that $\|\nabla\cQ u_j\|_{L^2}=\|\div u_j\|_{L^2},$ we have 
      \begin{multline}\label{eq:cLj}   \frac 12\frac d{dt}\cL_j^2+\mu\|\nabla\cP u_j\|_{L^2}^2+(1-\kappa)\|\nabla\cQ u_j\|_{L^2}^2
   +\kappa\|\nabla a_j\|_{L^2}^2 
    = \kappa\int\Delta\cQ u_j\cdot\nabla a_j\\
 +     \frac12\int(|u_j|^2-a_j^2)\div u+\kappa \int \div u \,\div(a_ju_j)
   -\kappa\int(u_j\cdot\nabla u)\cdot\nabla a_j\\
    +\int R_j^1 a_j+\int(R_j^2+g_j)\cdot u_j
   +\kappa\int u_j\cdot\nabla R_j^1   +\kappa\int (R_j^2+g_j)\cdot\nabla a_j .   \end{multline}
   Bernstein and Cauchy-Schwarz inequalities ensure that 
   $$   \biggl|\int\Delta\cQ u_j\cdot\nabla a_j\biggr|\leq 2^{j+1}\|\nabla \cQ u_j\|_{L^2}\|\nabla a_j\|_{L^2}.$$
Hence, for all $j_0\in\N,$ there exist $\kappa>0$ and  $c>0$ depending only on $j_0$ such that for all $j\leq j_0,$ we have
 \begin{eqnarray}\label{equiv:Lj} &\displaystyle\frac12\|(a_j,u_j)\|_{L^2} \leq \cL_j\leq 2 \|(a_j,u_j)\|_{L^2}\andf\\\nonumber
 & \displaystyle\mu\|\nabla\cP u_j\|_{L^2}^2\!+\!(1-\kappa)\|\nabla\cQ u_j\|_{L^2}^2
   \!+\!\kappa\|\nabla a_j\|_{L^2}^2 - \kappa\!\int\!\Delta\cQ u_j\!\cdot\!\nabla a_j\geq c\min(1,\mu)\|(\nabla a_j,\nabla u_j)\|_{L^2}^2.
    \end{eqnarray}
   By frequency localization, we have 
   $$   \|(\nabla a_j,\nabla u_j)\|_{L^2}\simeq 2^{j}\cL_j.$$
   Reverting to \eqref{eq:cLj} and integrating over $[0,t],$ this implies (changing $c$ if necessary) that
    \begin{equation}\label{eq:cLj2}   \cL_j(t)+c\min(1,\mu)2^{2j}\int_0^t\cL_j\leq \cL_j(0)
 +    C\biggl(\int_0^t\|\nabla u\|_{L^\infty}\cL_j + \int_0^t\|(R_j^1,R_j^2)\|_{L^2}+\int_0^t\|g_j\|_{L^2}\biggr)\cdotp   \end{equation}
The commutator terms $R_j^1$ and $R_j^2$ may be bounded thanks to Lemma 2.100 and Remark 2.102 of \cite{BCD} that give us
\begin{equation}\label{eq:Rj2}\sum_{j\in\Z} 2^{js} \|[u,\ddj]\cdot\nabla z\|_{L^2}\lesssim \|\nabla u\|_{\dot B^{\frac d2}_{2,1}}\|z\|_{\dot B^s_{2,1}}\quad\hbox{if }
-d/2<s\leq 1+d/2.\end{equation}
This allows us to bound $R_j^2$ directly.  For $R_j^1,$ we note that
$$R_j^1=[u,\ddj]\cdot\nabla a+ \ddj a\,\div u-\ddj(a\,\div u).$$ 
From \eqref{emb:crit}, the product law $\dot B^s_{2,1}\times \dot B^{\frac d2}_{2,1}\to\dot B^s_{2,1}$
for $-d/2<s\leq d/2,$  and \eqref{eq:Rj2}, we also get
$$\sum_{j\in\Z} 2^{js} \|R_j^1\|_{L^2}\lesssim \|\nabla u\|_{\dot B^{\frac d2}_{2,1}}\|a\|_{\dot B^s_{2,1}}.$$
Reverting to \eqref{eq:cLj2}, multiplying by $2^{js},$  summing
 up on $j\leq j_0+1$ and observing that \eqref{equiv:Lj} implies\footnote{For the time being,
 the low and high frequencies parts of the norms correspond to \eqref{eq:not1} and \eqref{eq:not2} with 
 $\alpha=2^{j_0}$ for  some sufficiently large positive integer  $j_0$. We will explain later how to get the same result with $j_0=0.$}
 \begin{equation}\label{eq:equiv1} \sum_{j\leq j_0} 2^{js}\cL_j\simeq \|(a,u)\|_{\dot B^s_{2,1}}^\ell,
\end{equation}
  we end up with
\begin{multline}\label{eq:lf1}\|(a,u)(t)\|_{\dot B^s_{2,1}}^\ell+\!\int_0^t \!\|(a,u)\|_{\dot B^{s+2}_{2,1}}^\ell
\lesssim \|(a_0,u_0)\|_{\dot B^s_{2,1}}^\ell +\! \int_0^t\!\|\nabla u\|_{\dot B^{\frac d2}_{2,1}} \|(a,u)\|_{\dot B^s_{2,1}}
+\!\int_0^t\!\|g\|_{\dot B^s_{2,1}}^\ell.\end{multline}
As the function $Q$ defined in \eqref{eq:Q} satisfies $Q(0)=0$ and $Q'(0)=1,$  it is of the form $Q(a)=a+ak(a),$ and  thus
(see \eqref{eq:parabolic}):
\begin{equation}\label{eq:decompow}
w=w^h+u^\ell+(-\Delta)^{-1}\nabla a^\ell +(-\Delta)^{-1}\nabla(a k(a))^\ell,
\end{equation}
whence
\begin{align*}g
&=k_1(a)\cA w^h+k_1(a)\cA u^\ell-k_1(a)\nabla a^\ell  -k_1(a)\nabla (a k_2(a))^\ell  
+\nabla (a k_3(a)). \end{align*}
In what follows, we shall use repeatedly that  whenever $s\leq s',$ we have
       \begin{equation}\label{eq:lfhf} \|z \|_{\dot B^{s'}_{2,1}}^{\ell,\alpha}\lesssim \alpha^{s'-s} \|z \|_{\dot B^{s}_{2,1}}^{\ell,\alpha}
\andf \|z\|^{h,\alpha}_{\dot B^{s}_{2,1}} \lesssim \alpha^{s-s'} 
 \|z\|^{h,\alpha}_{\dot B^{s'}_{2,1}}.\end{equation}
Using product laws in Besov spaces, and the stability of the space $\dot B^{\frac d2}_{2,1}$ by left-composition
(see \cite[Chap. 2]{BCD} and remember that $\|a\|_{L^\infty}$ is small), this implies that for all  $-d/2<s\leq d/2-1,$ 
\begin{align}\label{eq:aka}
\|\nabla(ak_3(a))\|_{\dot B^s_{2,1}}^\ell &\leq \|\nabla(a^\ell k_3(a))\|_{\dot B^s_{2,1}}^\ell+\|\nabla(a^hk_3(a))\|_{\dot B^s_{2,1}}^\ell
\nonumber\\
 &\lesssim  \|a^\ell k_3(a)\|_{\dot B^{s+1}_{2,1}}+2^{j_0}\|a^hk_3(a)\|_{\dot B^s_{2,1}}
 \nonumber\\&\lesssim \bigl(\|a\|^\ell_{\dot B^{s+1}_{2,1}}+\|a\|^h_{\dot B^s_{2,1}}\bigr)\|a\|_{\dot B^{\frac d2}_{2,1}}.\end{align}
Next, from product laws and composition estimates,  and \eqref{eq:aka}, we have
\begin{align}
\label{eq:cA1} \Bigl\|k_1(a) \cA w^h\|_{\dot B^s_{2,1}} &\lesssim \|a\|_{\dot B^{\frac d2}_{2,1}}
\|w\|^h_{\dot B^{s+2}_{2,1}},\\
 \label{eq:cA2}\Bigl\|k_1(a) \cA u^\ell\|_{\dot B^s_{2,1}} &\lesssim \|a\|_{\dot B^{\frac d2}_{2,1}}
\|u\|^\ell_{\dot B^{s+2}_{2,1}},\\ 
 \label{eq:cA3}
  \Bigl\|k_1(a)\nabla a^\ell\|_{\dot B^{s}_{2,1}}&\lesssim\|a\|_{\dot B^{\frac d2}_{2,1}}\|\nabla a\|^\ell_{\dot B^s_{2,1}} ,
   \\\label{eq:cA4}
   \Bigl\|k_1(a)\nabla(ak_2(a))^\ell\|_{\dot B^{s}_{2,1}}&\lesssim\|a\|_{\dot B^{\frac d2}_{2,1}}
   \|\nabla(ak_2(a))\|_{\dot B^s_{2,1}}^\ell \lesssim  \|a\|_{\dot B^{\frac d2}_{2,1}}^2\bigl(\|a\|^\ell_{\dot B^{s+1}_{2,1}}+\|a\|^h_{\dot B^s_{2,1}}\bigr)\cdotp\end{align}
  Hence,  reverting to \eqref{eq:lf1} and remembering \eqref{eq:smalla}, we conclude that
  \begin{multline}\label{eq:lf2}\|(a,u)(t)\|_{\dot B^s_{2,1}}^\ell+\!\int_0^t \!\|(a,u)\|_{\dot B^{s+2}_{2,1}}^\ell
\lesssim \|(a_0,u_0)\|_{\dot B^s_{2,1}}^\ell +\! \int_0^t\!\|\nabla u\|_{\dot B^{\frac d2}_{2,1}} \|(a,u)\|_{\dot B^s_{2,1}}
\\+\!\int_0^t\!\|a\|_{\dot B^{\frac d2}_{2,1}}\bigl(\|u\|_{\dot B^{s+2}_{2,1}}^\ell+\|w\|_{\dot B^{s+2}_{2,1}}^h\bigr)
+\int_0^t\|a\|_{\dot B^{\frac d2}_{2,1}}\bigl(\|a\|^\ell_{\dot B^{s+1}_{2,1}}+\|a\|^h_{\dot B^s_{2,1}}\bigr)\cdotp\end{multline}

      \subsubsection*{Estimates for the  high frequencies}
      
      To get optimal estimates for the high frequencies of the solution, it is  wise
      to rewrite \eqref{eq:systemlf} in terms of $a$ and $w$ (defined in \eqref{eq:parabolic}). Now, because  $\div u=\div w+Q(a),$   
      the mass equation may be written
      \begin{equation}\label{eq:aa}   \d_ta +u\cdot \nabla a+F(a)=-\div w-a\,\div w\with F(z):=(1+z)Q(z).\end{equation}
Hence, observing that   $F(0)=0$ and $F'(0)=1,$ we have 
$$         \d_ta +u\cdot \nabla a+a=-\div w-a\,\div u +a k(a).$$
Applying $\ddj$ to the above equation gives
$$  \d_ta_j +u\cdot \nabla a_j+a_j= [u,\ddj]\cdot\nabla a-\div w_j- \ddj(a\,\div u +a k(a)).$$
Therefore, we easily get 
\begin{multline*}\|a_j(t)\|_{L^2}+\int_0^t\|a_j\|_{L^2}\leq \|a_j(0)\|_{L^2}+\int_0^t\|\div u\|_{L^\infty}\|a_j\|_{L^2}\\
+\int_0^t\bigl(\| [u,\ddj]\cdot\nabla a\|_{L^2}
+\|\div w_j\|_{L^2}+\|\ddj(a\div u +a k(a))\|_{L^2}\bigr)\cdotp\end{multline*}
The commutator term may be bounded in light of \eqref{eq:Rj2} and we use \eqref{eq:lfhf} to handle the next term.
Hence multiplying by $2^{js}$ and summing on $j\geq j_0$ yields 
\begin{multline*}
\|a(t)\|^h_{\dot B^s_{2,1}}+\int_0^t \|a\|^h_{\dot B^s_{2,1}}\leq \|a_0\|^h_{\dot B^s_{2,1}}+\int_0^t\|\div u\|_{L^\infty}\|a\|^h_{\dot B^s_{2,1}}
\\+
C2^{-j_0}\int_0^t\|w\|_{\dot B^{s+2}_{2,1}}^h+
C\int_0^t\|\nabla u\|_{\dot B^{\frac d2}_{2,1}}\|a\|_{\dot B^s_{2,1}}+
\int_0^t\|a\,\div u+ ak(a)\|_{\dot B^s_{2,1}}^h.\end{multline*}
The term corresponding to $a\,\div u$ may be bounded exactly as the commutator. Furthermore,  we observe that, 
owing to the high frequency cut-off and to \eqref{eq:lfhf}, 
\begin{align}\label{eq:aka1}
\| ak(a)\|_{\dot B^s_{2,1}}^h&\leq \| a^\ell k(a)\|_{\dot B^s_{2,1}}^h+\|a^hk(a)\|_{\dot B^s_{2,1}}^h\nonumber\\\nonumber
&\leq  2^{-j_0}\| a^\ell k(a)\|_{\dot B^{s+1}_{2,1}}+\|a^hk(a)\|_{\dot B^s_{2,1}}\\
&\lesssim \bigl(\|a\|^\ell_{\dot B^{s+1}_{2,1}}+\|a\|^h_{\dot B^s_{2,1}}\bigr)\|a\|_{\dot B^{\frac d2}_{2,1}}.
\end{align}
Therefore, we have
\begin{multline}\label{eq:a}
\|a(t)\|^h_{\dot B^s_{2,1}}+\int_0^t \|a\|^h_{\dot B^s_{2,1}}\leq \|a_0\|^h_{\dot B^s_{2,1}}\\+
C2^{-j_0}\int_0^t\|w\|_{\dot B^{s+2}_{2,1}}^h+
C\int_0^t\|\nabla u\|_{\dot B^{\frac d2}_{2,1}}\|a\|_{\dot B^s_{2,1}}+
C\int_0^t  \bigl(\|a\|^\ell_{\dot B^{s+1}_{2,1}}+\|a\|^h_{\dot B^s_{2,1}}\bigr)\|a\|_{\dot B^{\frac d2}_{2,1}}.\end{multline}
With regard to  $w,$ from the velocity and mass equations, we discover that
$$
\partial_tw+\cA w=-u\cdot\nabla u +\frac a{1+a}\cA w-(-\Delta)^{-1}\nabla\bigl(Q'(a)\div((1+a) u)\bigr)\cdotp$$
Since $Q'(0)=1$  and $\div w=\div u-Q(a),$ we get  with our  convention for functions $k_i$: 
\begin{equation}\label{eq:decompo}Q'(a)\div((1+a) u)=\div w+a +a k_1(a)+ k_2(a)\div u+(1+k_3(a))\div(au).\end{equation}
Hence, one can recast the equation for $w$ as:
\begin{multline*}
\partial_tw+\cA w+(-\Delta)^{-1}\nabla\bigl(\div w+a)=G\\
\with G:=-u\cdot\nabla u+ k_1(a)\cA w - (-\Delta)^{-1}\nabla\bigl(a k_2(a)+ k_3(a)\div u+(1+k_4(a))\div(au)\bigr)\cdotp
\end{multline*}
Applying once again $\ddj$ to the equation 
 and remembering \eqref{eq:redond},
we easily get 
\begin{equation*}\|w(t)\|_{\dot B^s_{2,1}}^h + \min(1,\mu)
\int_0^t\|w\|_{\dot B^{s+2}_{2,1}}^h \leq \|w_0\|_{\dot B^s_{2,1}}^h 
+ \int_0^t\bigl(\|w\|_{\dot B^{s}_{2,1}}^h+\|a\|_{\dot B^{s-1}_{2,1}}^h\bigr)
+\int_0^t\|G\|^h_{\dot B^s_{2,1}}
.\end{equation*}
To bound $G,$ we just have to use classical product and composition estimates, \eqref{eq:lfhf}  and  \eqref{eq:aka1},
and, as regards the term with $\cA w,$ Relation \eqref{eq:decompow}, then Inequalities \eqref{eq:cA1} to \eqref{eq:cA4}. We get
\begin{align*}
\|u\cdot\nabla u\|_{\dot B^s_{2,1}}^h&\lesssim \|u\|_{\dot B^s_{2,1}}\|\nabla u\|_{\dot B^{\frac d2}_{2,1}},\\
\Bigl\|k_1(a)\,\cA w\|_{\dot B^s_{2,1}}^h&\lesssim 
\|a\|_{\dot B^{\frac d2}_{2,1}}\bigl(\|u\|_{\dot B^{s+2}_{2,1}}^\ell+\|w\|_{\dot B^{s+2}_{2,1}}^h\bigr)\\&\hspace{1.6cm}+
\|a\|_{\dot B^{\frac d2}_{2,1}}\|a\|^\ell_{\dot B^{s+1}_{2,1}}
+\|a\|_{\dot B^{\frac d2}_{2,1}}^2
\bigl(\|a\|^\ell_{\dot B^{s+1}_{2,1}}+\|a\|^h_{\dot B^s_{2,1}}\bigr),\\
\|(-\Delta)^{-1}\nabla(ak_2(a))\|_{\dot B^s_{2,1}}^h&\lesssim 2^{-j_0}\|ak_2(a)\|_{\dot B^s_{2,1}}^h
\lesssim  \bigl(\|a\|^\ell_{\dot B^{s+1}_{2,1}}+\|a\|^h_{\dot B^s_{2,1}}\bigr)\|a\|_{\dot B^{\frac d2}_{2,1}},\\
\|(-\Delta)^{-1}\nabla(k_3(a)\div u)\|_{\dot B^s_{2,1}}^h&\lesssim 2^{-j_0}\|k_3(a)\div u\|_{\dot B^s_{2,1}}^h
\lesssim \|a\|_{\dot B^{\frac d2}_{2,1}}\|u\|_{\dot B^{s+1}_{2,1}},\\ 
\|(-\Delta)^{-1}\nabla((1+k_4(a))\div(au))\|_{\dot B^s_{2,1}}^h&\lesssim
2^{-j_0}\|(1+k_4(a))\div(au)\|_{\dot B^s_{2,1}}\\
&\lesssim 2^{-j_0}(1+\|a\|_{\dot B^{\frac d2}_{2,1}}) \|au\|_{\dot B^{s+1}_{2,1}}\\
&\lesssim 2^{-j_0}(1+\|a\|_{\dot B^{\frac d2}_{2,1}}) \|a\|_{\dot B^{\frac d2}_{2,1}} \|u\|_{\dot B^{s+1}_{2,1}}.
\end{align*}
Hence, using \eqref{eq:lfhf} to bound the linear terms, then \eqref{eq:smalla}  yields
\begin{multline}\label{eq:wh}\|w(t)\|_{\dot B^s_{2,1}}^h + 
\int_0^t\|w\|_{\dot B^{s+2}_{2,1}}^h \lesssim \|w_0\|_{\dot B^s_{2,1}}^h + 
 \int_0^t\bigl(2^{-2j_0}\|w\|_{\dot B^{s+2}_{2,1}}^h+2^{-j_0}\|a\|_{\dot B^{s}_{2,1}}^h\bigr)\\
+\int_0^t\|\nabla u\|_{\dot B^{\frac d2}_{2,1}}\|u\|_{\dot B^s_{2,1}}
+\int_0^t\|a\|_{\dot B^{\frac d2}_{2,1}}\bigl(\|u\|_{\dot B^{s+2}_{2,1}}^\ell+\|w\|_{\dot B^{s+2}_{2,1}}^h+
\|a\|^\ell_{\dot B^{s+1}_{2,1}}+\|a\|^h_{\dot B^s_{2,1}}+\|u\|_{\dot B^{s+1}_{2,1}}\bigr)\cdotp
\end{multline}
Putting this inequality together with \eqref{eq:a}, taking $j_0$ large enough to eliminate the
linear terms in the right-hand side and remembering \eqref{eq:smalla}, we end up with 
\begin{multline}\label{eq:aw1}
\|a(t)\|^h_{\dot B^s_{2,1}}+\|w(t)\|_{\dot B^s_{2,1}}^h + 
\int_0^t\bigl(\|a\|^h_{\dot B^s_{2,1}}+\|w\|_{\dot B^{s+2}_{2,1}}^h\bigr)\lesssim 
 \|a_0\|^h_{\dot B^s_{2,1}}+ \|w_0\|_{\dot B^s_{2,1}}^h \\+ 
\int_0^t\|\nabla u\|_{\dot B^{\frac d2}_{2,1}}\|(a,u)\|_{\dot B^s_{2,1}}
+\int_0^t \|a\|_{\dot B^{\frac d2}_{2,1}}\bigl(\|u\|^\ell_{\dot B^{s+2}_{2,1}}+\|w\|^h_{\dot B^{s+2}_{2,1}}\bigr)\\
+\int_0^t\|a\|_{\dot B^{\frac d2}_{2,1}}\bigl(\|a\|_{\dot B^{s+1}_{2,1}}^\ell+ \|a\|_{\dot B^{s}_{2,1}}^h
+\|u\|_{\dot B^{s+1}_{2,1}}\bigr)\cdotp
\end{multline}
From \eqref{eq:decompow}, it is not difficult to see that
 $$\|u\|_{L^\infty_t(\dot B^{s}_{2,1})}^h\leq\|w\|_{L^\infty_t(\dot B^{s}_{2,1})}^h
 +C2^{-j_0} \bigl(\|a\|_{L^\infty_t(\dot B^{s}_{2,1})}^h +  \|a\|_{L^\infty_t(\dot B^{\frac d2}_{2,1})} \|a\|_{L^\infty_t(\dot B^{s}_{2,1})}\bigr)\cdotp$$
Hence,  since in addition to \eqref{eq:smalla}, the function $\nabla u$ is small in $L^1(\R_+;\dot B^{\frac d2}_{2,1}),$ 
we discover (by interpolation and Young inequality)  that the second line of \eqref{eq:aw1} can be absorbed  by the sum of the left-hand sides of \eqref{eq:lf2} and of \eqref{eq:aw1}. 
Next, we see that for all $\eps>0,$ we have
\begin{align*}
\int_0^t\|a\|_{\dot B^{\frac d2}_{2,1}}\|a\|^\ell_{\dot B^{s+1}_{2,1}}
&\lesssim \eps^{-1}\int_0^t \|a\|_{\dot B^{\frac d2}_{2,1}}^2\|a\|^\ell_{\dot B^{s}_{2,1}}
+\eps \int_0^t \|a\|^\ell_{\dot B^{s+2}_{2,1}}\andf\\
\int_0^t\|a\|_{\dot B^{\frac d2}_{2,1}}\|a\|^h_{\dot B^{s}_{2,1}}
&\lesssim \eps^{-1}\int_0^t \|a\|_{\dot B^{\frac d2}_{2,1}}^2\|a\|^h_{\dot B^{s}_{2,1}}
+\eps \int_0^t \|a\|^h_{\dot B^{s}_{2,1}}.\end{align*}
Furthermore, since $\nabla u=\nabla w-(-\Delta)^{-1}\nabla^2 Q(a)$ with 
$Q(a)=a+ak(a),$ 
we have, due to \eqref{eq:aka1} and interpolation
\begin{align*}
\|u\|_{\dot B^{s+1}_{2,1}}&\leq \|u\|_{\dot B^{s+1}_{2,1}}^\ell
+\|w\|_{\dot B^{s+1}_{2,1}}^h+C\|Q(a)\|_{\dot B^{s}_{2,1}}^h\\
&\lesssim  \bigl(\|u\|_{\dot B^{s}_{2,1}}^\ell\|u\|_{\dot B^{s+2}_{2,1}}^\ell\bigr)^{1/2}
+\bigl(\|w\|_{\dot B^{s}_{2,1}}^h\|w\|_{\dot B^{s+2}_{2,1}}^h\bigr)^{1/2}+
\|a\|_{\dot B^{s}_{2,1}}^h+ \bigl(\|a\|^\ell_{\dot B^{s+1}_{2,1}}+\|a\|^h_{\dot B^s_{2,1}}\bigr)\|a\|_{\dot B^{\frac d2}_{2,1}}.
\end{align*}
Hence the norm of $u$ in $L^2(0,t;\dot B^{s+1}_{2,1})$ may be bounded by the sum of the left-hand sides of \eqref{eq:lf2} and of \eqref{eq:aw1}.
One can now conclude that there exists $\eta>0$ such that if
\begin{equation}\label{eq:smallcritic}\|a\|_{L^\infty(\dot B^{\frac d2}_{2,1})}+ \|a\|_{L^2(\dot B^{\frac d2}_{2,1})}+ 
\|\nabla u\|_{L^1(\dot B^{\frac d2}_{2,1})}\leq\eta,\end{equation}
 then we have for all $-d/2<s\leq d/2-1$ and $t\in\R_+,$ 
\begin{multline}\label{eq:aw2}
\|(a,u)\|_{L_t^\infty(\dot B^s_{2,1})}+\|w\|^h_{L_t^\infty(\dot B^s_{2,1})}+ \|(a,u)\|^\ell_{L^1_t(\dot B^{s+2}_{2,1})}+\|u\|_{L^2_t(\dot B^{s+1}_{2,1})}\\+
 \|a\|^h_{L^1_t(\dot B^{s}_{2,1})}+\|w\|^h_{L^1_t(\dot B^{s+2}_{2,1})}
\lesssim  \|(a_0,u_0)\|_{\dot B^s_{2,1}}.
\end{multline}
Although our previous computations lead to a threshold between low and high frequencies 
located at $2^{j_0}$ for some large enough positive integer $j_0,$ 
it can be shifted at $j_0=0.$  In fact, we only have to check that the blocks $\ddj w$ for $j\in\{0,\cdots,j_0\}$
are bounded in $L^1(0,t;L^2)$ and $L^\infty(0,t;L^2)$ independently of $t.$
This can be done from   \eqref{eq:parabolic} and the fact that $Q(a)=a+ak(a),$  which imply that
$$\ddj w=\ddj u+(-\Delta)^{-1}\nabla\ddj a+(-\Delta)^{-1}\nabla\ddj(ak(a)).$$
For $j\in\{0,\cdots,j_0-1\},$ the first two terms may be bounded thanks to \eqref{eq:aw2}. 
For the last term, we note that $a$ is bounded in $(L^2\cap L^\infty)(\R_+;\dot B^{\frac d2}_{2,1}),$ hence
$ak(a),$ in $(L^1\cap L^\infty)(\R_+;\dot B^{\frac d2}_{2,1}).$

\subsubsection*{The borderline case $s=-d/2$} 
Propagating  regularity $\dot B^{-\frac d2}_{2,\infty}$ is needed
for justifying the uniform convergence of $a$ to $b$ in the case $d=2.$
The strategy is the same as before once we observe  that, although  the product law 
$\dot B^{\frac d2}_{2,1}\times \dot B^{s}_{2,1}$ as well as the commutator estimate \eqref{eq:Rj2} fail 
for $s=-d/2,$ one can consider  the Besov norm:
\begin{equation}\label{eq:besovinfini}
\|z\|_{\dot B^{-\frac d2}_{2,\infty}}:=\sup_{j\in\Z} 2^{-j\frac d2}\|\ddj z\|_{L^2} \end{equation} since we have
$$\| z_1\,z_2\|_{\dot B^{-\frac d2}_{2,\infty}} \lesssim \| z_1\|_{\dot B^{\frac d2}_{2,1}} 
\| z_2\|_{\dot B^{-\frac d2}_{2,\infty}} 
\andf \sup_{j\in\Z} 2^{-j\frac d2}\|[u,\ddj]\cdot\nabla z\|_{L^2}\lesssim \|\nabla u\|_{\dot B^{\frac d2}_{2,1}}
\| z\|_{\dot B^{-\frac d2}_{2,\infty}} .$$
Instead of summing on $j\leq j_0+1$ to get \eqref{eq:lf1}, 
we take the supremum, and similarly for the high frequencies. 
In the end, we obtain the following modification of Inequality \eqref{eq:aw2}: 
\begin{multline}\label{eq:awlimit}
\|(a,u)\|_{L_t^\infty(\dot B^{-\frac d2}_{2,\infty})}+\|w\|^h_{L_t^\infty(\dot B^{-\frac d2}_{2,\infty})}+ 
\|(a,u)\|^\ell_{\wt L^1_t(\dot B^{-\frac d2+2}_{2,\infty})}+\|u\|_{\wt L^2_t(\dot B^{-\frac d2+1}_{2,\infty})}\\+
 \|a\|^h_{\wt L^1_t(\dot B^{-\frac d2}_{2,\infty})}+\|w\|^h_{\wt L^1_t(\dot B^{-\frac d2+2}_{2,\infty})}
\lesssim  \|(a_0,u_0)\|_{\dot B^{-\frac d2}_{2,\infty}},
\end{multline}
where we used the following  notation for $\sigma\in\R$ and $1\leq q,r\leq\infty$:
\begin{equation}\label{eq:tilde}
\begin{array}{c}
\|z\|^\ell_{\wt L^q_t(\dot B^\sigma_{2,r})}:=\bigl\| 2^{js}\|\ddj z\|_{L^q(0,t;L^2)}\bigr\|_{\ell^r(j\leq 1)},
\quad \|z\|^h_{\wt L^q_t(\dot B^\sigma_{2,r})}:=\bigl\| 2^{js}\|\ddj z\|_{L^q(0,t;L^2)}\bigr\|_{\ell^r(j\geq0)}\\[2ex]
\andf   \|z\|_{\wt L^q_t(\dot B^\sigma_{2,r})}:=\bigl\| 2^{js}\|\ddj z\|_{L^q(0,t;L^2)}\bigr\|_{\ell^r}.\end{array}
\end{equation}
The details of the proof are left to the reader.

 \subsubsection*{Estimates in the diffusive scaling} 
 
In light of \eqref{eq:rescaling2} and \eqref{eq:rescalingb}, we  get  for all $t\in\R_+$ and $\sigma\in\R,$  
 $$ \|(a,u)(t)\|_{\dot B^\sigma_{2,1}}= \check \nu^{\sigma-\frac d2}  \Bigl\|\Bigl(\check a^\nu,\frac{\check u^\nu}{c\nu}\Bigr)\Bigl(\frac t{c^2}\Bigr)\Bigr\|_{\dot B^\sigma_{2,1}}
\!\! \andf\!\! \int_0^t\!\|(a,u)\|_{\dot B^{\sigma+2}_{2,1}}= \check \nu^{\sigma-\frac d2} \nu^2\int_0^{\frac t{c^2}}  
 \Bigl\|\Bigl(\check a^\nu,\frac{\check u^\nu}{c\nu}\Bigr)\Bigr\|_{\dot B^\sigma_{2,1}}.$$ 
 Hence, Inequality \eqref{eq:aw2} implies that, if 
 \begin{equation}\label{eq:smallcriticnu}\|\check a^\nu\|_{L^\infty(\dot B^{\frac d2}_{2,1})}+ c \|\check a^\nu\|_{L^2(\dot B^{\frac d2}_{2,1})}+ \|\nabla \check u^\nu\|_{L^1(\dot B^{\frac d2}_{2,1})}\leq\eta\ll1,\end{equation}
then  (remember that $u_0^\nu=\nu^{-1}\check u_0^\nu$ and $a_0^\nu=\check a_0^\nu$):
 \begin{multline}\label{eq:awnu}
c\|\check a^\nu\|_{L_t^\infty(\dot B^s_{2,1})}+ \nu^{-1}\Bigl(\|\check u^\nu\|_{L_t^\infty(\dot B^s_{2,1})}+\|\check w^\nu\|^{h,\check\nu^{-1}}_{L_t^\infty(\dot B^s_{2,1})}\Bigr)
+\|\check u\|_{L^2_t(\dot B^{s+1}_{2,1})}+\nu \|\check u^\nu\|^{\ell,\check\nu^{-1}}_{L^1_t(\dot B^{s+2}_{2,1})}\\+\nu \|\check w^\nu\|^{h,\check\nu^{-1}}_{L^1_t(\dot B^{s+2}_{2,1})}+
c\nu^2\|\check a^\nu\|^{\ell,\check\nu^{-1}}_{L^1_t(\dot B^{s+2}_{2,1})}+c^3\|\check a^\nu\|^{h,\check\nu^{-1}}_{L^1_t(\dot B^{s}_{2,1})}
\lesssim   Y^s_{\nu,c,0}:=c\|a_0^\nu\|_{\dot B^s_{2,1}}+ \|u_0^\nu\|_{\dot B^s_{2,1}}.
\end{multline}
Owing to \eqref{eq:uniform2}, 
Condition \eqref{eq:smallcriticnu}  is ensured if \eqref{eq:smalldata} holds true.
 Similarly, after rescaling, \eqref{eq:awlimit} gives
\begin{multline}\label{eq:awnulimit}
c\|\check a^\nu\|_{L_t^\infty(\dot B^{-\frac d2}_{2,\infty})}+ \nu^{-1}\Bigl(\|\check u^\nu\|_{L_t^\infty(\dot B^{-\frac d2}_{2,\infty})}
+\|\check w^\nu\|^{h,\check\nu^{-1}}_{L_t^\infty(\dot B^{-\frac d2}_{2,\infty})}\Bigr)
+\|\check u^\nu\|_{\wt L^2_t(\dot B^{-\frac d2+1}_{2,\infty})}
+\nu \|\check u^\nu\|^{\ell,\check\nu^{-1}}_{\wt L^1_t(\dot B^{-\frac d2+2}_{2,\infty})}\\+\nu \|\check w^\nu\|^{h,\check\nu^{-1}}_{\wt L^1_t(\dot B^{-\frac d2+2}_{2,\infty})}+
c\nu^2\|\check a^\nu\|^{\ell,\check\nu^{-1}}_{\wt L^1_t(\dot B^{-\frac d2+2}_{2,\infty})}
+c^3\|\check a^\nu\|^{h,\check\nu^{-1}}_{\wt L^1_t(\dot B^{-\frac d2}_{2,\infty})}
\lesssim   c\|a_0^\nu\|_{\dot B^{-\frac d2}_{2,\infty}}+ \|u_0^\nu\|_{\dot B^{-\frac d2}_{2,\infty}}.
\end{multline}


      \subsubsection*{Step 3: Study of the limit equation}
   Here we want to solve \eqref{eq:limit} supplemented with initial data $b_0$ in $\dot B^{\frac d2}_{2,1}\cap \dot B^s_{2,1}.$ 
   This equation may be rewritten 
      \begin{eqnarray}\label{eq:limitb}&\d_tb +v\cdot\nabla b+c^2 F(b)=0\\[1ex]  \label{eq:v}
      &\with    v:=-c^2(-\Delta)^{-1}\nabla Q(b)\andf      F(z):=(1+z) Q(z).\end{eqnarray}
     Note the similarity with the vorticity equation of two-dimensional incompressible perfect flows. 
     \smallbreak
     As a first, let us derive a priori estimates in Besov spaces for a (smooth) solution of \eqref{eq:limitb}.
     Decomposing $F(b)$ into $F(b)=b -bk(b)$ and localizing the above equation in the Fourier space
     by means of  $(\ddj)_{j\in\Z},$ we discover that $b_j:=\ddj b$ satisfies
     $$\d_tb_j+v\cdot\nabla b_j+c^2  b_j= c^2\ddj(b k(b))+[v,\ddj]\cdot\nabla b.$$
    Hence, implementing the usual  energy method on the above
    localized equation and observing that  $\div v=c^2 Q(b),$  we  get for all $j\in\Z$ and $t\in\R_+,$ 
    $$    \|b_j(t)\|_{L^2}+c^2\int_0^t\|b_j\|_{L^2}\leq \|b_j(0)\|_{L^2}
    +c^2\int_0^t\bigl(\|b_j\,Q(b)\|_{L^2}+\|\ddj(b k(b))\|_{L^2}\bigr)+\int_0^t\|[v,\ddj]\cdot\nabla b\|_{L^2}.$$
    Let   $\sigma\in(-d/2,d/2].$
     Multiplying the above inequality by $2^{j\sigma},$ summing on $j\in\Z$ and using the commutator estimate \eqref{eq:Rj2} and the definition of $v$  yields, 
     with the notation defined in \eqref{eq:tilde},    
     $$   \|b\|_{\wt L^\infty_t(\dot B^\sigma_{2,1})}+c^2\int_0^t\|b\|_{\dot B^\sigma_{2,1}}\leq \|b_0\|_{\dot B^\sigma_{2,1}}
     +\int_0^t\Bigl(c^2\bigl(\|Q(b)\|_{L^\infty}+C\|\nabla v\|_{\dot B^{\frac d2}_{2,1}}\bigr)\|b\|_{\dot B^\sigma_{2,1}}+c^2\|bk(b)\|_{\dot B^\sigma_{2,1}}\Bigr)\cdotp
     $$
     Finally using composition lemma, product law, the definition of $v$ in \eqref{eq:v}, 
       the  critical embedding \eqref{emb:crit} 
        and the fact that  $\nabla^2(-\Delta)^{-1}$ is a self-map on $\dot B^\sigma_{2,1}$
     allows us  to conclude that 
      $$   \|b\|_{\wt L^\infty_t(\dot B^{\sigma}_{2,1})}+c^2\int_0^t\|b\|_{\dot B^\sigma_{2,1}}\leq \|b_0\|_{\dot B^\sigma_{2,1}}
     +Cc^2\int_0^t\|b\|_{\dot B^{\frac d2}_{2,1}}\|b\|_{\dot B^\sigma_{2,1}}\quad\hbox{for all }\ \sigma\in[s,d/2].$$
     From this inequality and a standard iterative scheme, it is 
     easy to show that if $\|b_0\|_{\dot B^{\frac d2}_{2,1}}$ is small enough, 
     then \eqref{eq:limitb} admits a unique global solution $b\in\cC_b(\R_+;\dot B^{\frac d2}_{2,1})$
     and that the regularity $\dot B^s_{2,1}$ (for any $s\in(-d/2,d/2])$ is preserved. Furthermore, we have
      \begin{equation}\label{eq:bs}
\|b\|_{\wt L^\infty_t(\dot B^{\sigma}_{2,1})}+\frac{c^2}2\int_0^t\|b\|_{\dot B^{\sigma}_{2,1}}\leq \|b_0\|_{\dot B^{\sigma}_{2,1}}
\quad\hbox{for all }\ \sigma\in[s,d/2]. \end{equation}
      We claim that global well-posedness holds true whenever 
         \begin{equation}\label{eq:limitb0}  b_0\in \dot B^{\frac d2-1}_{2,1}\cap\dot B^{\frac d2}_{2,1}\with 1+b_0>0.\end{equation}
     Indeed, let  $T$ be the lifespan of the corresponding solution  $b$  seen as a continuous function of the
     time variable, with range in  $\dot B^{\frac d2-1}_{2,1}\cap\dot B^{\frac d2}_{2,1}.$
Then,  $v$ defined in \eqref{eq:v}
  is  in $\cC([0,T);\dot B^{\frac d2}_{2,1}\cap \dot B^{\frac d2+1}_{2,1}).$ 
 In light of the Cauchy-Lipschitz theorem and of \eqref{emb:crit},
     the vector-field $v$ has a unique $C^1$ flow $X_v$   on $[0,T)\times\R^d.$
    Furthermore,  $\wt b:=b\circ X_v$ satisfies the following ordinary differential equation: 
     \begin{equation}\label{eq:limitwtb}\d_t\wt b+c^2F(\wt b)=0.\end{equation}
   If  $Q(z)\geq 0$ for $z\geq0,$ and $Q(z)\leq 0$ for $z\in(-1,0),$ then 
     the $L^\infty$ norm  of $\wt b$  is nonincreasing on $[0,T),$ and $1+\wt b(t,\cdotp)\geq \inf_{x\in\R^d}(1+b_0(x)).$ 
    Hence, if  $\|\wt b(t)\|_{L^\infty}$ and $\inf_x \wt b(t,x)$  control  the blow-up of the $\dot B^{\sigma}_{2,1}$ norm
    for $\sigma\in\{-1+\frac d2,\frac d2\},$  then   \eqref{eq:limitb} is globally well-posed . 
 Now, since  \eqref{eq:limitwtb}  may be rewritten 
    $$    \d_t\wt b+c^2\wt b=c^2\wt b\, k(\wt b) \with k(0)=0,$$
    localizing  by means of $(\ddj)_{j\in\Z}$ and using an energy method yields
    $$    \|\wt b(t)\|_{\dot B^{\sigma}_{2,1}}+c^2\int_0^t \|\wt b\|_{\dot B^{\sigma}_{2,1}}\leq \|b_0\|_{\dot B^{\sigma}_{2,1}}+
    c^2\int_0^t\|\wt b k(\wt b)\|_{\dot B^{\sigma}_{2,1}}\quad\hbox{for all }\ \sigma\in[-d/2,d/2].$$
    Tame estimates for the product and composition in  the space $\dot B^{\frac d2}_{2,1}$ (see \cite{BCD}) ensure that 
    $$\|z \,k(z)\|_{\dot B^{\frac d2}_{2,1}}\leq C_{m,M} \|z\|_{L^\infty}\|z\|_{\dot B^{\frac d2}_{2,1}}\quad\hbox{if}\quad
    \|z\|_{L^\infty}\leq M<\infty\andf 1+z\geq m>0.$$
    Hence, taking advantage of Gronwall lemma, one can conclude that 
    if $\wt b$ (or, equivalently, $b$) is bounded on $[0,T)\times\R^d,$ then it may be continued beyond $T.$
    Finally, since, as emphasized before, $\|\wt b(t)\|_{L^\infty}$ is nonincreasing,  any solution to \eqref{eq:limitb} 
    supplemented with initial data $b_0$ satisfying \eqref{eq:limitb0} is  global, and belongs to $\cC(\R_+;\dot B^{\frac d2}_{2,1}\cap\dot B^{\frac d2-1}_{2,1} ).$
    Reverting to the original unknown $b$ allows to conclude that \eqref{eq:limitb} is globally well-posed for 
    data satisfying \eqref{eq:limitb0}, that any regularity $\dot B^\sigma_{2,1}$ with $\sigma\in[s,d/2]$ 
    is preserved, and that \eqref{eq:bs} holds true. 
       \subsubsection*{Step 4:  Passing to the limit by the direct Eulerian approach}
       In this step, we drop index $\nu$ on $\ca^\nu$ and $\cu^\nu$, for better 
       readability. Furthermore, we  just note $z^\ell$ or $z^h$ instead of $z^{\ell,\check\nu^{-1}}$
    and $z^{h,\check\nu^{-1}}$ (ditto for the norms).

Let $\da:=\check a-b.$ 
 Using   the limit velocity defined in \eqref{eq:v}, we see that the evolution of $\da$  is governed by the transport equation:
  \begin{equation}\label{eq:da1} \d_t\da
  + \check u\cdot\nabla\da +c^2(F(\check a)-F(b))=-(1+\check a)\div \check w +     (v-\check u)\cdot\nabla b.\end{equation}
   As $b$ has  regularity $\dot B^{\frac d2}_{2,1}$ and the last term causes a loss of one derivative, 
   we can only hope to control $\da(t)$ in $\dot B^{\frac d2-1}_{2,1}$ (this is the usual curse of derivative loss  in stability estimates for hyperbolic systems).
      Now, since $F'(0)=1,$ we have 
   $$F(\check a)-F(b)=\da(1-G(\check a,b))\with G(\check a,b):= (b-\check a)\int_0^1(1-\tau)F''(\check a+\tau \da)\,d\tau.$$ 
      Let us  use the following decompositions (with $k(0)=0$):
       \begin{align}\label{eq:dec1}
       \div \check w&=\div \check w^h+\div \check u^\ell -c^2\check a^\ell-c^2(\check a k(\check a))^\ell,\\\label{eq:dec2}
      v-\check u&=(v-\check u)^\ell-\check w^h+c^2 (-\Delta)^{-1}\nabla\bigl(Q(\check a)-Q(b)\bigr)^h\end{align}
       so as to rewrite \eqref{eq:da1} as: 
          \begin{align}\label{eq:da2} 
          \d_t\da+ \check u\cdot\nabla\da &+c^2\da=A+B \\
 \with         A&:=c^2(1+\check a)\check a^\ell\andf\nonumber\\ 
          B&:=c^2\da\, G(\check a,b) -(1+\check a)\div \check w^h-(1+\check a)\div \check u^\ell
+c^2(1+\check a)(\check ak(\check a))^\ell  \nonumber\\&\hspace{3cm}+\bigl((v-\check u)^\ell-\check w^h
+c^2 (-\Delta)^{-1}\nabla\bigl(Q(\check a)-Q(b)\bigr)^h\bigr)\cdot\nabla b.\nonumber\end{align}
        In the rest of the section, we assume that $s\in[d/2-2,d/2-1]$ (or  $s\in(-1,0]$ if $d=2$) and 
        define the Lebesgue exponents $p_s\in[1,2]$ and $q_s\in[2,\infty]$ from the relations
        \begin{equation}\label{eq:pq}     \frac1{p_s}=\frac d4-\frac s2\andf \frac1{q_s}=\frac1{p_s}-\frac12\cdotp\end{equation}
        Note that for  $s=d/2-2,$ we have $p_s=1$ and $q_s=2.$ 
        The definition of $p_s$ and $q_s$  is motivated by the fact that
        \begin{equation}\label{eq:lplq} \|z\|_{L_t^{p_s}(\dot B^{\frac d2}_{2,1})} \lesssim     \|z\|_{L_t^{1}(\dot B^{s+2}_{2,1})}^{\frac1{p_s}}
         \|z\|_{L_t^{\infty}(\dot B^{s}_{2,1})}^{\frac1{p'_s}}\andf
          \|z\|_{L_t^{q_s}(\dot B^{\frac d2-1}_{2,1})} \lesssim     \|z\|_{L_t^{1}(\dot B^{s+2}_{2,1})}^{\frac1{q_s}}
         \|z\|_{L_t^{\infty}(\dot B^{s}_{2,1})}^{\frac1{q'_s}}.\end{equation}
     Applying  to \eqref{eq:da1} the  estimates for the transport equation derived in e.g. \cite[Chap. 3]{BCD}, we get
        \begin{multline}\label{eq:da3}\|\da\|_{L_t^\infty(\dot B^{\frac d2-1}_{2,1})}
        +c^{\frac2{q_s}}\|\da\|_{L_t^{q_s}(\dot B^{\frac d2-1}_{2,1})}
\leq \|\da_0\|_{\dot B^{\frac d2-1}_{2,1}}+\int_0^t\|\nabla \check u\|_{\dot B^{\frac d2}_{2,1}}\|\da\|_{\dot B^{\frac d2-1}_{2,1}}\\+c^{-\frac2{q'_s}}\|A\|_{L^{q_s}_t(\dot B^{\frac d2-1}_{2,1})}
+c^{-\frac2{p'_s}}\|B\|_{L^{p_s}_t(\dot B^{\frac d2-1}_{2,1})}.\end{multline}
    To estimate $A,$ we use the fact that
    $$
    \|(1+\check a)\check a^\ell\|_{\dot B^{\frac d2-1}_{2,1}}\lesssim \bigl(1+\|\check a\|_{\dot B^{\frac d2}_{2,1}}\bigr)
    \|\check a\|^\ell_{\dot B^{\frac d2-1}_{2,1}}.$$
    Hence, taking advantage of  \eqref{eq:lplq}, of \eqref{eq:awnu} and of the smallness of $\|\check a\|_{L^\infty_t(\dot B^{\frac d2}_{2,1})}$ yields 
    \begin{align}\label{eq:A}   \|A\|_{L_t^{q_s}(\dot B^{\frac d2-1}_{2,1})}&\lesssim 
     c^2\bigl(\|\check a\|^\ell_{L_t^{1}(\dot B^{s+2}_{2,1})}\bigr)^{\frac1{q_s}}
         \bigl(\|\check a\|^\ell_{L_t^{\infty}(\dot B^{s}_{2,1})}\bigr)^{\frac1{q'_s}}\nonumber \\
         &\lesssim c^{1-\frac2{q'_s}} \nu^{-\alpha_s} Y^s_{\nu,c,0}\with \alpha_s:=\frac2{q_s}\cdotp
         \end{align}
         Let us next estimate $B$ in $L^{p_s}(\R_+;\dot B^{\frac d2-1}_{2,1}).$  To handle  the first term of $B,$ we write that
         $$
        c^{-\frac2{p'_s}}\: c^2 \|\da\, G(\check a,b)\|_{L_t^{p_s}(\dot B^{\frac d2-1}_{2,1})}\lesssim  
        \Bigl(c^{\frac2{q_s}}\|\da\|_{L^{q_s}_t(\dot B^{\frac d2-1}_{2,1})}\Bigr)
        \Bigl(c \|(\check a,b)\|_{L_t^2(\dot B^{\frac d2}_{2,1})}\Bigr)\cdotp$$
         Hence, as    $c\|(\check a,b)\|_{L^2_t(\dot B^{\frac d2}_{2,1})}$ is small owing to \eqref{eq:uniform2}, this term may be absorbed by the left-hand side of \eqref{eq:da2}.
       Next, using product and composition laws in Besov spaces, 
        Inequality \eqref{eq:awnu} and  \eqref{eq:lplq}, we can write that
       \begin{equation*}
       \|(1+\check a)\div(u^\ell+ \check w^h)\|_{L^{p_s}_t(\dot B^{\frac d2-1}_{2,1})}\lesssim 
        \bigl(1+\|\check a\|_{L^\infty_t(\dot B^{\frac d2}_{2,1})}\bigr)\bigl(\|\div \check u\|^\ell_{L^{p_s}_t(\dot B^{\frac d2-1}_{2,1})}
     +  \|\div \check w\|^h_{L^{p_s}_t(\dot B^{\frac d2-1}_{2,1})}\bigr)\cdotp  \end{equation*}
          From \eqref{eq:awnu} and  \eqref{eq:lplq}, we directly have 
       \begin{equation}\label{eq:bla}\|\check u\|^\ell_{L^{p_s}_t(\dot B^{\frac d2}_{2,1})}
       +\|\check w\|^h_{L^{p_s}_t(\dot B^{\frac d2}_{2,1})}
       \lesssim  \nu^{-\alpha_s}Y^s_{\nu,c,0}.\end{equation}
       Hence we conclude that 
       $$           \|(1+\check a)\div (u^\ell+\check w^h)\|_{L^{p_s}_t(\dot B^{\frac d2-1}_{2,1})}\lesssim 
            \nu^{-\alpha_s} Y^s_{\nu,c,0}.$$
       Next, still thanks to \eqref{eq:awnu} and  \eqref{eq:lplq},
              \begin{align*}
        c^2\|(1+\check a)(\check a^\ell k(\check a))^\ell\|_{L^{p_s}_t(\dot B^{\frac d2-1}_{2,1})}&\lesssim 
           c^2\bigl(1+\|\check a\|_{L^\infty_t(\dot B^{\frac d2}_{2,1})}\bigr)\|\check a\|^\ell_{L^{q_s}_t(\dot B^{\frac d2-1}_{2,1})}
           \|\check a\|_{L^2_t(\dot B^{\frac d2}_{2,1})}\\
           & \lesssim \nu^{-\alpha_s} \eta_0 Y^s_{\nu,c,0}.\end{align*}
        Using also \eqref{eq:lfhf}, we get 
 \begin{align*}
      c^2  \|(1+\check a)(\check a^h k(\check a))^\ell\|_{L^{p_s}_t(\dot B^{\frac d2-1}_{2,1})}
   &\lesssim c^2\bigl(1+\|\check a\|_{L^\infty_t(\dot B^{\frac d2}_{2,1})}\bigr)\|\check a^h k(\check a)\|^\ell_{L^{p_s}_t(\dot B^{\frac d2-1}_{2,1})}\\
   &\lesssim     c^2\check\nu^{-\alpha_s} \|\check a^h k(\check a)\|_{L^{p_s}_t(\dot B^{s}_{2,1})}\\
                    &\lesssim   c^{2+\frac2{q_s}}\nu^{-\alpha_s}   
                    \|\check a\|^h_{L^{q_s}_t(\dot B^{s}_{2,1})} \|\check a\|_{L^2_t(\dot B^{\frac d2}_{2,1})}
      \leq   \nu^{-\alpha_s} \eta_0 Y^s_{\nu,c,0}.\end{align*}
        Keeping  \eqref{eq:lplq}  and \eqref{eq:bs} in mind, 
         if we assume that
       \begin{equation}\label{eq:smallb0}
       \|b_0\|_{\dot B^{\frac d2}_{2,1}}\leq\eta_0\ll1\end{equation}  
       then  we get 
                         $$       \|\check w^h\cdot\nabla b\|_{L^{p_s}_t(\dot B^{\frac d2-1}_{2,1})}\lesssim
                 \|\check w^h\|_{L^{p_s}_t(\dot B^{\frac d2}_{2,1})}\|\nabla b\|_{L^{\infty}_t(\dot B^{\frac d2-1}_{2,1})}
                 \lesssim \nu^{-\alpha_s} Y^s_{\nu,c,0}\|b_0\|_{\dot B^{\frac d2}_{2,1}}\lesssim \nu^{-\alpha_s} \eta_0 Y^s_{\nu,c,0}.        $$
                             $$
       \|(v-\check u)^\ell\cdot\nabla b\|_{L^{p_s}_t(\dot B^{\frac d2-1}_{2,1})}\lesssim 
       \bigl(\|\check u\|^\ell_{L^{p_s}_t(\dot B^{\frac d2}_{2,1})}+\|v\|^\ell_{L^{p_s}_t(\dot B^{\frac d2}_{2,1})}\bigr)
       \|\nabla b\|_{L^{\infty}_t(\dot B^{\frac d2-1}_{2,1})}.$$
     Observe that \eqref{eq:lfhf} and \eqref{eq:bs} give us
       $$       \|v\|^\ell_{L^{p_s}_t(\dot B^{\frac d2}_{2,1})}\lesssim c^2 \check\nu^{-\alpha_s}
              \|Q(b)\|_{L^{p_s}_t(\dot B^{s}_{2,1})}\lesssim c\nu^{-\alpha_s}\|b_0\|_{\dot B^s_{2,1}}.$$
       Hence, using also \eqref{eq:bla}
        $$\|(v-\check u)^\ell\cdot\nabla b\|_{L^{p_s}_t(\dot B^{\frac d2-1}_{2,1})}\lesssim 
      \nu^{-\alpha_s}\bigl(Y^s_{\nu,c,0}+c\|b_0\|_{\dot B^s_{2,1}}\bigr) \|b_0\|_{\dot B^{\frac d2}_{2,1}}
      \leq   \nu^{-\alpha_s}\bigl(Y^s_{\nu,c,0}+c\|b_0\|_{\dot B^s_{2,1}}\bigr) \eta_0.$$
       Finally, we have 
                 \begin{align*}
       \| \bigl((-\Delta)^{-1}\nabla\bigl(Q(\check a)-Q(b)\bigr)^h\bigr)\cdot\nabla b\|_{L^{p_s}_t(\dot B^{\frac d2-1}_{2,1})}&\!\lesssim 
       \| (-\Delta)^{-1}\nabla\bigl(Q(\check a)-Q(b)\bigr)\|^h_{L^{\infty}_t(\dot B^{\frac d2}_{2,1})}
       \|\nabla b\|_{L^{p_s}_t(\dot B^{\frac d2-1}_{2,1})}\\
       &\!\lesssim \|\da\|_{L^{\infty}_t(\dot B^{\frac d2-1}_{2,1})}\|b\|_{L^{p_s}_t(\dot B^{\frac d2}_{2,1})} \\
          &\!\lesssim c^{-\frac2{p_s}} \|\da\|_{L^{\infty}_t(\dot B^{\frac d2-1}_{2,1})}\|b_0\|_{\dot B^{\frac d2}_{2,1}}.
       \end{align*}
    Hence, thanks to \eqref{eq:smallb0},  this latter term may be absorbed by the left-hand side of \eqref{eq:da2}.
    Finally, putting together \eqref{eq:A} and the above estimates pertaining to $B$ in \eqref{eq:da3} leads to
            \begin{equation}\label{eq:da5}\|\da\|_{L_t^\infty(\dot B^{\frac d2-1}_{2,1})}
        +c^{\frac2{q_s}}\|\da\|_{L_t^{q_s}(\dot B^{\frac d2-1}_{2,1})}
\leq \|\da_0\|_{\dot B^{\frac d2-1}_{2,1}}+ C\check\nu^{-\alpha_s} \bigl(c^{-1}Y^s_{\nu,c,0}+\|b_0\|_{\dot B^s_{2,1}}\bigr)\cdotp\end{equation}
  In order to prove the uniform convergence of $a$ to $b,$ we need to control additionally the high frequencies  of $\da$ in $L^{p_s}(\R_+;\dot B^{\frac d2-1}_{2,1})$.
                  For technical reasons, we are going to shift slightly on the right  the threshold
                  between low and high frequencies. This modification will be indicated using the exponents $\ell'$ and  $h'$ rather than $\ell$ and $h$.
                  \smallbreak
                  Now, applying   $\ddj$ to \eqref{eq:da2} then the usual energy argument and the commutator estimate \eqref{eq:Rj2}, 
          then summing  on $j>1+\log_2\check\nu,$  we discover that
          \begin{equation}\label{eq:dashift} c^{\frac2{p_s}}\|\da\|^{h'}_{L_t^{p_s}(\dot B^{\frac d2-1}_{2,1})}
\lesssim \|\da_0\|_{\dot B^{\frac d2-1}_{2,1}}^{h'}
+c^{-\frac2{p'_s}}\|A+B\|_{L^{p_s}_t(\dot B^{\frac d2-1}_{2,1})}^{h'}
+\int_0^t\|\nabla\check u\|_{\dot B^{\frac d2}_{2,1}}\|\da\|_{\dot B^{\frac d2-1}_{2,1}}.\end{equation}
                   The term  $B$ is estimated as above. For $A,$ we note that
          $$A=c^2 \check a^\ell + c^2\check a \check a^\ell.$$
Since $\ddj a^\ell=0$ for  $j>1+\log_2\check\nu,$ 
the linear   part of $A$ has no contribution to \eqref{eq:dashift}. As for the second part, it may be just estimated from 
                    \begin{equation}\label{eq:star0} c^2\|\check a \check a^\ell\|_{L^{p_s}_t(\dot B^{\frac d2-1}_{2,1})}\lesssim   
          c\bigl(c\|\check a\|_{L^2_t(\dot B^{\frac d2}_{2,1})}\bigr)  \|\check a\|^\ell_{L^{q_s}_t(\dot B^{\frac d2-1}_{2,1})}\lesssim 
          \nu^{-\alpha_s}\eta_0 Y^s_{\nu,c,0}.\end{equation}
          Hence, we have, due to \eqref{eq:dashift}, \eqref{eq:star0}  and  the previous bounds for $B,$
          \begin{equation}\label{eq:da6}
                  c^{\frac2{p_s}}\|\da\|^{h'}_{L_t^{p_s}(\dot B^{\frac d2-1}_{2,1})}
\lesssim \|\da_0\|_{\dot B^{\frac d2-1}_{2,1}}^{h'} +  \check\nu^{-\alpha_s} \bigl(c^{-1}Y^s_{\nu,c,0}+\|b_0\|_{\dot B^s_{2,1}}\bigr)\cdotp
\end{equation}
This inequality entails that $u$ converges to $v$ in $L^{p_s}(\R_+;\dot B^{\frac d2}_{2,1}).$ Indeed, 
consider  the decomposition \eqref{eq:dec2} \emph{with the threshold between 
low and high frequencies shifted on the right} (as explained above). Then, in light of  \eqref{eq:awnu} and \eqref{eq:lplq},  we have
\begin{equation}\label{eq:adfs}
\|\check u\|_{L_t^{p_s}(\dot B^{\frac d2}_{2,1})}^{\ell'} + \|v\|_{L_t^{p_s}(\dot B^{\frac d2}_{2,1})}^{\ell'} 
+\|\check w\|_{L_t^{p_s}(\dot B^{\frac d2}_{2,1})}^{h'}\lesssim \nu^{-\alpha_s}\bigl(Y_{\nu,c,0}^s+c\|b_0\|_{\dot B^s_{2,1}}\bigr)\cdotp\end{equation}
The result  is obvious for $\check w$ since 
$\|\check w\|_{L_t^{p_s}(\dot B^{\frac d2}_{2,1})}^{h'}\leq \|\check w\|_{L_t^{p_s}(\dot B^{\frac d2}_{2,1})}^{h}$.
It is also clear for $v,$ owing to \eqref{eq:bs} with $s=d/2-1.$ 
For $\check u,$ Inequalities   \eqref{eq:awnu} and \eqref{eq:lplq} allow to control 
$\|\check u\|_{L_t^{p_s}(\dot B^{\frac d2}_{2,1})}^{\ell}.$ 
To handle the dyadic bloc corresponding to  $j_{\check\nu}:=\lfloor \log_2\check\nu^{-1}\rfloor,$
we write that
$$\begin{aligned}
\|\dot\Delta_{j_{\check\nu}} \check u\|_{L_t^{p_s}(\dot B^{\frac d2}_{2,1})}
&\lesssim \|\dot\Delta_{j_{\check\nu}} \check u\|_{L_t^{1}(\dot B^{s+2}_{2,1})}^\theta
\|\dot\Delta_{j_{\check\nu}} \check u\|_{L_t^{\infty}(\dot B^s_{2,1})}^{1-\theta}
\with \theta:=\frac d4-\frac s2\cdotp\end{aligned}$$
Since one can shift the threshold in \eqref{eq:awnu} on the right (up to a change of the constant $C$),  one can conclude that 
\eqref{eq:adfs}  holds.

It is now easy to establish that $\check u$ converges to $v$ in $L^{p_s}(\R_+;\dot B^{\frac d2}_{2,1}).$
The starting point is  \eqref{eq:dec2} shifted on the right (that is, $\ell$ and $h$ are replaced by $\ell'$ and $h',$ respectively).
The first two terms  are bounded by means of \eqref{eq:adfs}. 
 For the last one, we write that
$Q(b)-Q(\check a)=\da(1+H(\check a,b))$ for some smooth function $H$ vanishing at $(0,0).$ Hence
\begin{align}\nonumber\|(-\Delta)^{-1}\bigl(\nabla(Q(\check a)-Q(b))\bigr)^{h'}\|_{L^{p_s}_t(\dot B^{\frac d2}_{2,1})}
&\lesssim \|Q(\check a)-Q(b)\|^{h'}_{L^{p_s}_t(\dot B^{\frac d2-1}_{2,1})}\\
&\lesssim \|\da\|^{h'}_{L^{p_s}_t(\dot B^{\frac d2-1}_{2,1})}+
\|\da\|_{L^{q_s}_t(\dot B^{\frac d2-1}_{2,1})}\|(\check a,b)\|_{L^2_t(\dot B^{\frac d2}_{2,1})}.\label{eq:star} \end{align}
Hence, taking advantage of \eqref{eq:da5} and \eqref{eq:da6}, we conclude that
\begin{equation}\label{eq:u-v}
\|\check u-v\|_{L_t^{p_s}(\dot B^{\frac d2}_{2,1})}\lesssim 
c^{\frac2{q'_s}}\|\da_0\|_{\dot B^{\frac d2-1}_{2,1}}+\nu^{-\alpha_s}\bigl(Y_{\nu,c,0}^s+c\|b_0\|_{\dot B^s_{2,1}}\bigr)\cdotp
\end{equation}
At this stage, one can specify the convergence of the solution of $(CNS_{\nu,c})$ in the original scaling as
we have 
$$
a(t)-b(\nu^{-1}t) =\da(\nu^{-1} t)\andf     u(t)-\nu^{-1}v(\nu^{-1}t)=\nu^{-1}(\check u-v)(\nu^{-1} t).$$
From \eqref{eq:da5}, we readily deduce \eqref{eq:asympa} 
while \eqref{eq:u-v} implies  \eqref{eq:asympv}. 

  \subsubsection*{Step 5:  The Lagrangian approach, and the uniform convergence}
              
In this step, we prove that if $a_0\to b_0$ in $\dot B^{\frac d2}_{2,1}$ and $(a_0^\nu,u_0^\nu)$ is bounded
in $\dot B^{\frac d2-2}_{2,1}$ (case $d\geq3$)  or  in $\dot B^{-1}_{2,\infty}$ (case $d=2$),  
 then the convergence of 
 $\check a$ to $b$ as $\nu$ goes to $\infty$  is  \emph{uniform}  on $\R_+\times\R^d.$
   %
   The key to the proof is  to  consider \eqref{eq:rescaledmass} in the Lagrangian coordinates
   pertaining  to the flow of the vector-field $\check u$ and to look at the convergence toward  $\wt b$ solution of \eqref{eq:limitwtb}.
   In this way, we can avoid
   the loss of one derivative, and  get  convergence in $\dot B^{\frac d2}_{2,1}$ but in Lagrangian coordinates.
   \medbreak
   More precisely, 
    let $X_{\check u}$   be the $C^1$ flow of the vector-field $\check u$ and let  
     $\wt a:=\check a \circ X_{\check u}.$ Then, $\wt a$ satisfies 
    \begin{equation}\label{eq:limitwta}\d_t\wt a + c^2 F(\wt a)=-(1+\wt a)\bigl((\div \check w)\circ X_{\check u}\bigr).\end{equation}
    It is nothing more than  a damped ordinary differential equation
    with a  source term  which, presumably, will tend to zero
    when $\nu$ goes to infinity. 
          Now,   $\wt\da:=\wt a-\wt b$ satisfies:
     $$ \d_t\wt\da+c^2\bigl(F(\wt a)-F(\wt b)\bigr)=-(1+\check a\circ X_{\check u})\bigl((\div \check w)\circ X_{\check u}\bigr)\cdotp$$
     Using the fact that $F'(0)=0$ and the mean value formula
$$F(z_1)-F(z_2) =(z_1-z_2) F'(z_2+\phi (z_1-z_2)),\qquad \phi \in(0,1), $$ 
reveals that if  $\ca$ and $b$ are small enough  in $\dot B^{\frac d2}_{2,1}$  and if 
$$-(1+\check a\circ X_{\check u})\bigl((\div \check w)\circ X_{\check u}\bigr)=\wt A+\wt B\with \wt A\in L^2(\R_+;\dot B^{\frac d2}_{2,1})\andf \wt B\in L^1(\R_+;\dot B^{\frac d2}_{2,1}),$$
then we have
\begin{equation}\label{eq:da7}\|\wt\da\|_{L^\infty_t(\dot B^{\frac d2}_{2,1})}+\frac c2 \|\wt\da\|_{L^2_t(\dot B^{\frac d2}_{2,1})}
  \leq  \|\da_0\|_{\dot B^{\frac d2}_{2,1}}+
     \|\wt B\|_{L^1_t(\dot B^{\frac d2}_{2,1})}
  +c^{-1}  \|\wt A\|_{L^2_t(\dot B^{\frac d2}_{2,1})}.\end{equation}
  Keeping \eqref{eq:dec1} in mind, we set 
  $$
  \wt A=-c^2\bigl((1+\check a) \check a^\ell\bigr)\circ X_{\check u}\andf
  \wt B=\bigl((1+\check a)(\div\check w^h+\div\check u^\ell-c^2(\check a k(\check a))^\ell)\bigr)\circ X_{\check u}.$$
  Using \cite{D14} (for the stability of $\dot B^{\frac d2}_{2,1}$ by right-composition) and product laws, we discover that
  $$
  \|\wt A\|_{\dot B^{\frac d2}_{2,1}}\lesssim c^2(1+\|\check a\|_{\dot B^{\frac d2}_{2,1}})\|\check a^\ell\|_{\dot B^{\frac d2}_{2,1}}.$$
  Hence, thanks to  \eqref{eq:smalla} and to \eqref{eq:awnu} with $s=d/2-1,$ we get
  \begin{equation}\label{eq:da8}
  c^{-1}  \|\wt A\|_{L^2_t(\dot B^{\frac d2}_{2,1})}\lesssim c\|\check a^\ell\|_{L^2_t(\dot B^{\frac d2}_{2,1})}\lesssim
  \nu^{-1} Y^{\frac d2-1}_{\nu,c,0}.\end{equation}
  In order to bound $\wt B$ we use also \eqref{eq:smalla} and the stability of $\dot B^{\frac d2}_{2,1}$ by right-composition, to get 
    $$  \|\wt B\|_{\dot B^{\frac d2}_{2,1}}\lesssim 
  \|\div \check u\|^{\ell,\check\nu^{-1}}_{\dot B^{\frac d2}_{2,1}} 
+ \|\div \check w\|^{h,\check\nu^{-1}}_{\dot B^{\frac d2}_{2,1}} 
  +c^2\|(\check a k(\check a))^\ell\|_{\dot B^{\frac d2}_{2,1}}.$$
  From \eqref{eq:awnu} with $s=d/2-1,$ we readily have
$$\|\div \check u\|^{\ell,\check\nu^{-1}}_{L_t^1(\dot B^{\frac d2}_{2,1})} 
+ \|\div \check w\|^{h,\check\nu^{-1}}_{L_t^1(\dot B^{\frac d2}_{2,1})} 
\lesssim \nu^{-1}Y_{\nu,c,0}^{\frac d2-1}.$$
From the same inequality, it is easy to deduce that 
$$\|\check a\|^\ell_{L^2_t(\dot B^{\frac d2}_{2,1})} \lesssim c^{-1}\nu^{-1} Y_{\nu,c,0}^{\frac d2-1}\andf
\|\check a\|^{h}_{L^2_t(\dot B^{\frac d2-1}_{2,1})} \lesssim  c^{-2}Y_{\nu,c,0}^{\frac d2-1}$$
while \eqref{eq:uniform2} implies that
$$\|\check a\|^\ell_{L^2_t(\dot B^{\frac d2}_{2,1})} \lesssim c^{-1}\eta_0.$$
Hence, using $a=a^\ell+a^h,$ then  composition and product laws and, finally,  \eqref{eq:lfhf}, we discover that
$$c^2\|(\check a k(\check a))^\ell\|_{L^2_t(\dot B^{\frac d2}_{2,1})}\lesssim \eta_0 \nu^{-1}Y_{\nu,c,0}^{\frac d2-1}$$
and thus
  \begin{equation*} \|\wt B\|_{L^1_t(\dot B^{\frac d2}_{2,1})}\lesssim
  \nu^{-1}Y_{\nu,c,0}^{\frac d2-1}.\end{equation*}
So,  combining with \eqref{eq:da7} and \eqref{eq:da8}, we  end up  for all $t\in\R_+$ with
\begin{equation}\label{eq:da9}\|\wt\da(t)\|_{\dot B^{\frac d2}_{2,1}}+c\|\wt\da\|_{L^2_t(\dot B^{\frac d2}_{2,1})}\lesssim \|\da_0\|_{\dot B^{\frac d2}_{2,1}}+
 \nu^{-1}Y^{\frac d2-1}_{\nu,c,0}.\end{equation}
 The above inequality combined with the critical embedding \eqref{emb:crit} 
 ensures that $\wt a:= \check a\circ X_{\check u}$ converges uniformly to $\wt b:=b\circ X_v$ on  $\R_+\times\R^d$
 whenever  $a_0^\nu\to b_0$ in $\dot B^{\frac d2}_{2,1}.$ 
 In order to show  that it entails that  $\check a$ converges uniformly to  $b,$  
 it is enough to prove 
  that  $\check a\circ X_{\check u}-b\circ X_{\cu}$ converges uniformly to $0.$
 
Now,  in light of the triangle inequality and of the definition of $\wt\da,$ we have
\begin{equation}\label{eq:da10}
\| \check a\circ X_{\check u}-b\circ X_{\check u}\|_{L^\infty}\leq \|\wt\da\|_{L^\infty} 
+\|b\circ X_v-b\circ X_{\check u}\|_{L^\infty}\end{equation}
If   $a_0\to b_0$ in $\dot B^{\frac d2}_{2,1}$ then $\wt\da$ tends to $0$ uniformly on $\R_+\times\R^d,$
owing to \eqref{eq:da9} and \eqref{emb:crit}. 
Proving that the second term of \eqref{eq:da10} also tends uniformly to $0$   relies on the facts 
that $b$ is uniformly continuous on $\R_+\times\R^d$ and that  $X_{\check u}$ converges to $X_v$ uniformly on $\R_+\times\R^d.$
For the first fact, we note that, as a consequence of $b\in \wt L^\infty(\R_+;\dot B^{\frac d2}_{2,1})$ 
and thus  of $b\in \wt L^\infty(\R_+;\dot B^{0}_{\infty,1}),$ the series $\sum \ddj b$ converges normally to $b$ on $\R_+\times\R^d.$
Furthermore, each term $\ddj b$ is uniformly continuous on $\R_+\times\R^d$ since we have from the equation \eqref{eq:limitb}
and product laws
$$\partial_t b\in L^1(\R_+;\dot B^{\frac d2-1}_{2,1})\hookrightarrow L^1(\R_+;\dot B^{-1}_{\infty,1}).$$
To prove the convergence of the flow, we observe that, by definition of $X_{\check u}$ and of $X_v,$ we have  for all $t\in\R_+$ and $y\in\R^d,$
$$X_v(t,y)-X_{\check u}(t,y)=\int_0^t \bigl(v(\tau,X_v(\tau,y))-\check u(\tau,X_v(\tau,y))\bigr)
+\int_0^t \bigl(\check u(\tau,X_v(\tau,y))-\check u(t,X_{\check u}(\tau,y))\bigr)\cdotp$$
Hence
$$\bigl|X_v(t,y)-X_{\check u}(t,y)\bigr|\leq\int_0^t\|v-\check u\|_{L^\infty}
+\int_0^t\|\nabla\check u\|_{L^\infty}\bigl|X_v(\tau,y)-X_{\check u}(\tau,y)\bigr|$$
 and Gronwall lemma gives us for all $t\geq0,$
 $$ \|X_v(t)-X_{\check u}(t)\|_{L^\infty} \leq \int_0^te^{\int_\tau^t\|\nabla \check u\|_{L^\infty}}\|v-\check u\|_{L^\infty}.$$
 In the case $d\geq3,$ one can propagate the regularity $s=d/2-2$ by the flow of the compressible 
 Navier-Stokes equations. Then, the exponent $p_s$ defined in \eqref{eq:pq} is equal to $1,$ 
  and \eqref{eq:u-v} combined with \eqref{emb:crit} ensure that $\check u-v$ converges to zero 
 in $L^1(\R_+;L^\infty),$   whence the uniform convergence of $X_{\check u}$ to $X_v,$ and thus of $\check a$ to $b.$
 \medbreak
 In the two-dimensional case, we are not allowed to take $s=-1$ and to deduce  directly from  \eqref{eq:u-v} 
 that $\check u-v\to 0$ in $L^1(\R_+;L^\infty).$
 However, one can take advantage of the following logarithmic interpolation inequality (see e.g. \cite[Chap. 10]{BCD}):
 $$ \|\check u-v\|_{L^1_t(L^\infty)}\lesssim  \|\check u-v\|_{L^1_t(\dot B^1_{2,1})}\lesssim  \|\check u-v\|_{\wt L^1_t(\dot B^1_{2,\infty})}
 \log\biggl(e+\frac{ \|\check u-v\|_{L^1_t(\dot B^2_{2,1})}}{\|\check u-v\|_{\wt L^1_t(\dot B^1_{2,\infty})}}\biggr)\cdotp$$
 Due to \eqref{eq:uniform2} and to \eqref{eq:v}--\eqref{eq:bs} with $\sigma=d/2,$ the
 numerator inside the logarithm is uniformly bounded with respect to $t$ and $\nu.$
 Hence, it suffices to establish that $\check u-v$ tends to $0$ in $\wt L^1(\R_+;\dot B^1_{2,\infty})$
 to conclude as in the higher dimensional case. 
  Again, we start from \eqref{eq:dec2}.  Inequality \eqref{eq:awnulimit} gives us
 $$
 \|\check u\|^{\ell'}_{\wt L^1_t(\dot B^1_{2,\infty})}+  \|\check w\|^{h'}_{\wt L^1_t(\dot B^1_{2,\infty})}\lesssim \nu^{-1}\bigl(c\|a_0\|_{\dot B^{-1}_{2,\infty}}
 +\|u_0\|_{\dot B^{-1}_{2,\infty}}\bigr)\cdotp$$
 Next,   adapting  \eqref{eq:lfhf} and \eqref{eq:bs} yields 
 $$\begin{aligned}  \|v\|^{\ell'}_{\wt L^1_t(\dot B^1_{2,\infty})}&\lesssim c^2\|Q(b)\|^{\ell'}_{\wt L^1_t(\dot B^0_{2,\infty})}\\
 &\lesssim c^2\check\nu^{-1}\|Q(b)\|^{\ell'}_{\wt L^1_t(\dot B^{-1}_{2,\infty})}\\
 &\lesssim   c^2\check\nu^{-1}\|b\|_{\wt L^1_t(\dot B^{-1}_{2,\infty})}\lesssim \nu^{-1}c\|b_0\|_{\dot B^{-1}_{2,\infty}}.\end{aligned}$$
 Finally, we may write 
  $$\begin{aligned}  \|(-\Delta)^{-1}\nabla(Q(\check a)-Q(b))^{h'}\|_{\wt L^1_t(\dot B^1_{2,\infty})}&\lesssim
  \|Q(\check a)-Q(b)\|^{h'}_{\wt L^1_t(\dot B^0_{2,\infty})}\\
  &\lesssim \|\da\|^{h'}_{\wt L^1_t(\dot B^0_{2,\infty})}+\|\da\|_{\wt L^2_t(\dot B^0_{2,\infty})}\|(\check a,b)\|_{L^2_t(\dot B^1_{2,1})}.\end{aligned}$$
 Hence, the problem reduces to proving the convergence to $0$ of $ \|\da\|^{h'}_{\wt L^1_t(\dot B^0_{2,\infty})}$
 and $ \|\da\|_{\wt L^2_t(\dot B^0_{2,\infty})}.$
 This is just a matter of checking that $A$  (resp. $B$)  defined in \eqref{eq:da2} are  of order $\check\nu^{-1}$ 
  in  $\wt L^2(\R_+;\dot B^0_{2,\infty})$  (resp.   $\nu^{-1}$ in $\wt L^1(\R_+;\dot B^0_{2,\infty})$), 
  which does not present any new difficulty:
   follow the calculation leading to \eqref{eq:da5}, using the endpoint product laws in Besov spaces, whenever it is needed.

    
    \section{The case of variable viscosity coefficients}\label{s:variable}
    
    This section aims at extending  the convergence result to the 
    physically relevant case where the viscosity coefficients depend smoothly on the density. 
    The governing equations then read
      $$   \left\{\begin{aligned}
    &\d_t \rho+\div( \rho u)=0,\\
   &\rho(\d_tu+u\cdot\nabla u) +\nu  \cA_\rho u+c^2\nabla(P(\rho))=0\end{aligned}\right.\leqno(\wt{CNS}_{\nu,c}) $$
    where     \begin{equation}\label{eq:cArho}
    \cA_\rho: u\mapsto -\div(2\mu(\rho)D(u))-\nabla(\lambda(\rho)\div u)\with (D(u))_{ij}:=\frac12\Bigl(\d_i u^j+\d_ju^i\Bigr)\cdotp
    \end{equation}
    As before,  the (given) pressure function $P$ is smooth and  such that $P(1)=0$ and $P'(1)=1.$
    We suppose that  the viscosity coefficients $\lambda$ and $\mu$ depend
    smoothly on $\rho>0$ and, in order for  the second order operator $\cA_\rho$ to be invertible, 
    that     \begin{equation}\label{eq:elliptic}
    \mu(\rho)>0\andf \lambda(\rho)+2\mu(\rho)>0.
    \end{equation}
       In accordance with the case studied before (which corresponds to $\mu(\rho)=\mu$ and $\lambda(\rho)=1-2\mu$),
        we make the normalization hypothesis $\lambda(1)+2\mu(1)=1$
   \medbreak
    Denoting $\rho=1+a$ and $Q(a):=P(1+a),$  the above equations in the case $c=\nu=1$ rewrite 
        \begin{equation}\label{eq:CNSv}\left\{\begin{aligned}
    &\d_t a+\div( (1+a)u)=0,\\
   &(1+ a)(\d_t u+ u\cdot\nabla u) +\cA_{1+a}  u+P'(1+a)\nabla  a=0.\end{aligned}\right. \end{equation}
    As in the constant viscosity case, 
    it  is known  (see e.g. \cite{DX17}) that if the initial data satisfy
    $$\|(a_0,u_0)\|_{\dot B^{\frac d2-1}_{2,1}}+\|a_0\|_{\dot B^{\frac d2}_{2,1}}\leq\eta_0,$$
    then \eqref{eq:CNSv}  has a unique global solution $(a,u)$ satisfying \eqref{eq:globcns}, and 
    using  the rescaling \eqref{eq:rescaling} allows to extend this result to general   $\nu,c>0.$ 
    \medbreak
Under  the diffusive rescaling defined in \eqref{eq:diffusive}, the density equation rewrites
    $$\d_t\check\rho+\cu\cdot\nabla \check\rho -c^2\check\rho\div\bigl(\cA_{\check\rho}^{-1}(\nabla(P(\check\rho)))\bigr)=-\check\rho\,\div \check w\with
    \cw:= \cu+c^2\cA_{\crho}^{-1}\nabla (P(\crho))    $$
while we have
$$\nu^{-2}\crho
\bigl(\d_t\check u+\check u\cdot\nabla\cu\bigr)+\cA_{\crho}\cw=0.$$
Hence, we may expect $\cw\to0$ and thus that $\cu\to v$ and $\crho\to q$ with 
$v=-c^2\cA_q^{-1}(\nabla(P(q))),$ and $q,$ a solution to the following nonlocal transport equation:
\begin{equation}\label{eq:q}\d_tq+v\cdot\nabla q-c^2q\,\div\cA_q^{-1}(\nabla(P(q))=0.\end{equation}
The rest of this section is devoted to justifying rigorously this heuristics. 
As the general strategy is the same as in the constant viscosity case, we just sketch the proof, highlighting the main 
differences.

\subsection{Inversion of the elliptic operator}

The first step is to check the inversibility of $\cA_\rho.$ 
We claim that   under hypothesis \eqref{eq:elliptic} with $(\rho-1)\in\dot B^{\frac d2}_{2,1},$
 whenever one takes some scalar function $m$ in $\dot B^\sigma_{2,1}$ for some $\sigma\in(-d/2,d/2]$
then, there exists a unique (up to constant)   vector-field $z$ such that  $\nabla z$ is in $\dot B^\sigma_{2,1}$ and
\begin{equation}\label{eq:carho}
\cA_\rho z=\nabla m.\end{equation}
Although we expect the result to be true  for all $\rho>0$ such that $a:=\rho-1$ belongs to $\dot B^{\frac d2}_{2,1},$
we assume here that  $a$ is small in $\dot B^{\frac d2}_{2,1},$  since only the small 
fluctuation case will be needed in what follows. 
Then, we claim that one can find a  mapping $\Phi: m\mapsto\nabla z$ 
with $z$ satisfying \eqref{eq:carho} and
\begin{equation}\label{eq:z} \| \Phi(m)\|_{\dot B^{\sigma}_{2,1}}\lesssim \|m\|_{\dot B^\sigma_{2,1}},\qquad \sigma\in(-d/2,d/2].\end{equation}
The basic idea is to rewrite \eqref{eq:carho} as 
$$
\cA_1z=\bigl(\cA_1-\cA_\rho\bigr)z+\nabla m.$$
We observe that operator $\cA_1$ restricted to potential vector-fields is nothing other than $-\Delta.$
Hence, if the solution exists, it has to satisfy
$$\nabla z=\nabla\cA_1^{-1}\Bigl(\div\bigl(2(\mu(\rho)-\mu(1))D(z)\bigl)+\nabla\bigl((\lambda(\rho)-\lambda(1))\div z\bigr)\Bigr)+\nabla  (-\Delta)^{-1} \nabla m.$$
The first term of the right-hand side is a combination of Riesz operators
acting on $k(a)\nabla z.$ Hence, for  $\sigma\in(-d/2,d/2],$ we have
$$
\|\nabla z\|_{\dot B^{\sigma}_{2,1}}\lesssim \|k(a)\nabla z\|_{\dot B^\sigma_{2,1}}+  \|m\|_{\dot B^\sigma_{2,1}}
\lesssim \|a\|_{\dot B^{\frac d2}_{2,1}} \|\nabla z\|_{\dot B^{\sigma}_{2,1}}+  \|m\|_{\dot B^\sigma_{2,1}}.
$$
Leveraging the smallness of $a$ gives  Inequality \eqref{eq:z}. 
\smallbreak
The  proof of existence of a map $\Phi: \dot B^{\sigma}_{2,1}\to \dot B^{\sigma}_{2,1}$ satisfying \eqref{eq:z}
stems from the fixed point theorem (that $\sigma\leq d/2$ ensures the completeness of $\dot B^{\sigma}_{2,1}$).
In fact, our problem is equivalent to finding a fixed point for
$$
 M\mapsto \nabla \cA_1^{-1} \Bigl(\div\bigl(2(\mu(\rho)-\mu(1))\Bigl(\textstyle{\frac{M+M^\top}2}\Bigr)\Bigr)+\nabla\bigl((\lambda(\rho)-\lambda(1)){\rm Tr} M\bigr)\Bigr)
 +\nabla (-\Delta)^{-1}\nabla m$$
 on the space of matrices $M$ with coefficients in $\dot B^s_{2,1}$  and such that $\d_kM_{ij}=\d_iM_{kj}$ for all
 $1\leq i,j,k\leq d.$ 
 \medbreak 
In what follows, $\Phi$ will be denoted by $\nabla\cA_\rho^{-1}.$ So, in particular, the effective velocity reads
\begin{equation}\label{eq:parabolicv}
w=u+\cA_\rho^{-1}\nabla(Q(a)).\end{equation}

\subsection{Control of lower order norms} Here we want to bound   $\|(a,u)(t)\|_{\dot B^s_{2,1}}$ 
for all $t>0$  and  $-d/2<s\leq d/2-1.$  
Thanks to \eqref{eq:rescaling2}, it suffices to consider the case $\nu=c=1,$
and thus the following equations:
\begin{equation}\label{eq:CNSvariable}
\left\{\begin{aligned}&\d_ta+\div((1+a)u)=0,\\
&(1+a)(\d_tu+u\cdot\nabla u)+\cA_{\rho}u+\nabla(Q(a))=0.\end{aligned}\right.\end{equation}

\subsubsection*{Control of the low frequencies} 
For the time being, we fix some frequency threshold $2^{j_0}$ with $j_0\in\N$ large enough, 
and the decomposition into low and high frequencies has to be understood in this meaning. 
Getting the same estimates with  $j_0=0$ can be achieved afterward, as in the constant viscosity case. 
\smallbreak
In terms of $w,$ the second equation of \eqref{eq:CNSvariable} may be rewritten
$$\d_tu+u\cdot\nabla u+\cA_{\rho} w=\frac a{1+a}\cA_\rho w.$$
Hence, observing (with the usual convention) that
$$\cA_\rho w=\cA_1u+\nabla a+(\cA_\rho-\cA_1)u+\nabla(a k(a)),$$
we discover that \eqref{eq:CNSvariable} may be rewritten  as \eqref{eq:systemlf} with, now, 
$$g=k_1(a)\cA_\rho w+(\cA_1-\cA_\rho)u +\nabla(ak(a))\cdotp$$
The last term of the right-hand side has already been treated in the constant viscosity case. 
The second one is of the form $\nabla(k(a)\nabla u).$ To handle it, we use the decompositions 
\begin{equation}\label{eq:decompou}
u=u^\ell+w^h - \bigl(\cA_\rho^{-1}\nabla(Q(a))\bigr)^h
\end{equation}
and 
\begin{equation}\label{eq:decompoQ}
Q(a)=a^\ell+a^h +ak(a).\end{equation}
Hence
\begin{multline*}
(\cA_1-\cA_\rho)u \simeq \nabla\bigl(k(a)\nabla(u^\ell+w^h)\bigr)
+\nabla\bigl(k(a)\bigl(\nabla\cA_\rho^{-1}\nabla a^\ell\bigr)^h\bigr)\\
+\nabla\bigl(k(a)\bigl(\nabla\cA_\rho^{-1}\nabla a^h\bigr)^h\bigr)
+\nabla\bigl(k(a)\bigl(\nabla\cA_\rho^{-1}\nabla(a k(a))\bigr)^h\bigr)\cdotp\end{multline*}
Since $s+1\leq d/2,$ we readily have 
\begin{equation}\label{eq:uw}\|\nabla\bigl(k(a)\nabla(u^\ell+w^h)\bigr)\|_{\dot B^{s}_{2,1}}\lesssim 
\|a\|_{\dot B^{\frac d2}_{2,1}}\|u^\ell+w^h\|_{\dot B^{s+2}_{2,1}}.\end{equation}
Next, since $\nabla\cA_\rho^{-1}\nabla$ is a self-map on $\dot B^s_{2,1}$ and on $\dot B^{s+1}_{2,1}$ (see \eqref{eq:z}),  
we have
\begin{align*}
\|\nabla\bigl(k(a)\bigl(\nabla\cA_\rho^{-1}\nabla a^\ell\bigr)^h\bigr)\|_{\dot B^s_{2,1}}
&\lesssim \|k(a)\bigl(\nabla\cA_\rho^{-1}\nabla a^\ell\bigr)^h\|_{\dot B^{s+1}_{2,1}}\lesssim\|a\|_{\dot B^{\frac d2}_{2,1}}\|a^\ell\|_{\dot B^{s+1}_{2,1}},\\
\|\nabla\bigl(k(a)\bigl(\nabla\cA_\rho^{-1}\nabla a^h\bigr)^h\bigr)\|_{\dot B^s_{2,1}}^\ell
&\lesssim \|k(a)\bigl(\nabla\cA_\rho^{-1}\nabla a^h\bigr)^h\|_{\dot B^s_{2,1}}
\lesssim \|a\|_{\dot B^{\frac d2}_{2,1}}\|a^h\|_{\dot B^{s}_{2,1}},\\
\|\nabla\bigl(k(a)\bigl(\nabla\cA_\rho^{-1}\nabla(a k(a))\bigr)^h\bigr\|_{\dot B^s_{2,1}}^\ell
&\lesssim \|k(a)\bigl(\nabla\cA_\rho^{-1}\nabla(a k(a))\bigr)^h\bigr\|_{\dot B^s_{2,1}}\lesssim 
 \|a\|_{\dot B^{\frac d2}_{2,1}}^2\|a\|_{\dot B^{s}_{2,1}}.\end{align*}
 Finally, using  again \eqref{eq:decompoQ}, and that
 $$ \cA_\rho w=\cA_\rho(u^\ell +w^h) + \cA_\rho\bigl(\cA_\rho^{-1}\nabla(Q(a))\bigr)^\ell,$$
 we see that 
 $$k_1(a)\cA_\rho w=k_1(a)\cA_\rho(u^\ell+w^h) 
 + k_1(a)\cA_\rho\bigl(\cA_\rho^{-1}\nabla a^\ell\bigr)^\ell
  + k_1(a)\cA_\rho\bigl(\cA_\rho^{-1}\nabla (a^h+a k(a))\bigr)^\ell.$$
 We have
 \begin{align*}
 \|k_1(a)\cA_\rho(u^\ell+w^h)\|_{\dot B^s_{2,1}}&\lesssim \|a\|_{\dot B^{\frac d2}_{2,1}}\|u^\ell+w^h\|_{\dot B^{s+2}_{2,1}},\\ 
 \|k_1(a)\cA_\rho\bigl(\cA_\rho^{-1}\nabla a^\ell\bigr)^\ell\|_{\dot B^s_{2,1}}&\lesssim \|a\|_{\dot B^{\frac d2}_{2,1}}\|a^\ell\|_{\dot B^{s+1}_{2,1}},\\
  \| k_1(a)\cA_\rho\bigl(\cA_\rho^{-1}\nabla (a^h+ak(a))\bigr)^\ell\|_{\dot B^s_{2,1}}&\lesssim \|a\|_{\dot B^{\frac d2}_{2,1}}
  \|\nabla\cA_\rho^{-1}\nabla (a^h+ak(a))\|_{\dot B^{s+1}_{2,1}}^\ell\\
  &\lesssim \|a\|_{\dot B^{\frac d2}_{2,1}}  \|\nabla\cA_\rho^{-1}\nabla (a^h+ak(a))\|_{\dot B^{s}_{2,1}}\\
    &\lesssim \|a\|_{\dot B^{\frac d2}_{2,1}}  \bigl(\|a\|_{\dot B^{s}_{2,1}}^h 
    + \|a\|_{\dot B^{\frac d2}_{2,1}} \|a\|_{\dot B^{s}_{2,1}}\bigr)\cdotp
\end{align*}
Hence
\begin{equation}\label{eq:Arhol}
\|k_1(a)\cA_\rho w\|_{\dot B^s_{2,1}}\lesssim  \|a\|_{\dot B^{\frac d2}_{2,1}}  
\bigl(\|u^\ell+w^h\|_{\dot B^{s+2}_{2,1}}+\|a\|_{\dot B^{s+1}_{2,1}}^\ell+ \|a\|_{\dot B^{s}_{2,1}}^h 
    + \|a\|_{\dot B^{\frac d2}_{2,1}} \|a\|_{\dot B^{s}_{2,1}}\bigr)\cdotp\end{equation}

Remembering  \eqref{eq:smalla},  one can  conclude that $g$ may be estimated by the right-hand side of 
\eqref{eq:lf2}.

\subsubsection*{Control of the high frequencies} 

Let $F(a):=(1+a)Q'(a).$ 
In terms of $a$ and $w,$ the mass equation rewrites
\begin{equation}\label{eq:abis}\d_t a+u\cdot \nabla a+(1+a)\div w+F(a)= G:=(1+a)\div\bigl((\cA_\rho^{-1}-\cA_1^{-1})\nabla(Q(a))\bigr).\end{equation}
Compared to \eqref{eq:aa}, only the term $G$ is new. 
Using  the short notation $Q$ for $Q(a),$ we have
$$\begin{aligned}
(\cA_\rho^{-1}-\cA_1^{-1})\nabla Q&=\cA_1^{-1}(\cA_1-\cA_\rho)\cA_\rho^{-1}\nabla Q\\
&=\cA_1^{-1}\Bigl(2\div\bigl((\mu(\rho)-\mu(1))D\bigl(\cA^{-1}_\rho\nabla Q\bigr)\bigr)
+\nabla\bigl((\lambda(\rho)-\lambda(1)\bigr)\div\bigl(\cA^{-1}_\rho\nabla Q\bigr)\bigr)\Bigr)\cdotp\end{aligned}$$
The right-hand side is a homogeneous multiplier of degree $-1$ acting on 
functions of type $k(a)\nabla\cA_\rho^{-1}\nabla Q$ with $k(0)=0.$  
From the decomposition $k(a)=k'(0)(a^\ell+a^h)+ ak_1(a)$ with $k_1(0)=0,$ we thus infer that
$$G\simeq (1+a)\cR\bigl(a^\ell \nabla \cA_\rho^{-1} \nabla Q\bigr)+ (1+a)\cR\bigl((a^h+ak_1(a))\nabla \cA_\rho^{-1} \nabla Q\bigr),$$
where $\cR$ stands for some combination of Riesz operators. 
Hence, we have
$$
\|G\|_{\dot B^s_{2,1}}^h\lesssim 
2^{-j_0}\|(1+a)\cR\bigl(a^\ell \nabla \cA_\rho^{-1} \nabla Q\bigr)\|_{\dot B^{s+1}_{2,1}}^h
+ \|(1+a)\cR\bigl((a^h+ak_1(a))\nabla \cA_\rho^{-1} \nabla Q\bigr)\|_{\dot B^s_{2,1}}^h.$$
Using  \eqref{eq:smalla}, \eqref{eq:z},  product and composition laws in Besov spaces and the definition of $Q,$  we deduce that
$$\begin{aligned}\|G\|_{\dot B^s_{2,1}}^h&\lesssim \|a^\ell\nabla\cA_\rho^{-1}\nabla Q\|_{\dot B^{s+1}_{2,1}}+\|a^h\nabla\cA_\rho^{-1}\nabla Q\|_{\dot B^s_{2,1}}+\|a k_1(a)\nabla\cA_\rho^{-1}\nabla Q\|_{\dot B^s_{2,1}}\\
&\lesssim\bigl(\|a\|_{\dot B^{s+1}_{2,1}}^\ell+ \|a\|_{\dot B^{s}_{2,1}}^h
+ \|a k_1(a)\|_{\dot B^{s}_{2,1}}\bigr)\|Q\|_{\dot B^{\frac d2}_{2,1}}\\
& \lesssim \bigl(\|a\|_{\dot B^{s+1}_{2,1}}^\ell+\|a\|_{\dot B^s_{2,1}}^h
+\|a\|_{\dot B^{\frac d2}_{2,1}}\|a\|_{\dot B^s_{2,1}}\bigr) \|a\|_{\dot B^{\frac d2}_{2,1}}.\end{aligned}$$
In the end,  \eqref{eq:a} becomes
\begin{multline*}
\|a(t)\|^h_{\dot B^s_{2,1}}+\int_0^t \|a\|^h_{\dot B^s_{2,1}}\leq \|a_0\|^h_{\dot B^s_{2,1}}+
C2^{-j_0}\int_0^t\|w\|_{\dot B^{s+2}_{2,1}}^h\\+
C\int_0^t\Bigl(\|\nabla u\|_{\dot B^{\frac d2}_{2,1}}\|a\|_{\dot B^s_{2,1}}+
 \bigl(\|a\|^\ell_{\dot B^{s+1}_{2,1}}+\|a\|^h_{\dot B^s_{2,1}}\bigr)\|a\|_{\dot B^{\frac d2}_{2,1}}
+\|a\|^2_{\dot B^{\frac d2}_{2,1}}  \|a\|_{\dot B^s_{2,1}}\Bigr)\cdotp\end{multline*}
In the variable viscosity case, the equation satisfied by $w$  reads
\begin{equation}\label{eq:wv}\d_t w+\cA_1 w=-u\cdot\nabla u+\frac a\rho\cA_\rho w+(\cA_1-\cA_\rho)w+\partial_t(\cA_\rho^{-1}\nabla Q(a)).\end{equation}
We have to bound the high frequencies of the right-hand side in $L^1(\R_+;\dot B^s_{2,1}).$
The first term has already been  treated in the constant viscosity case, and the second one can be bounded
according to \eqref{eq:Arhol}.
For the third one, we use the decomposition
$$
(\cA_1-\cA_\rho)w=(\cA_1-\cA_\rho)(u^\ell+w^h)+(\cA_1-\cA_\rho)\bigl(\cA_\rho^{-1}\nabla(Q(a))\bigr)^\ell.$$
The first part may be bounded according to \eqref{eq:uw}.
The last term is of the form $$\nabla(k(a)(\nabla\cA_\rho^{-1}\nabla(Q(a)))^\ell).$$ It can be bounded by decomposing $Q$ as in 
\eqref{eq:decompoQ}. Now, thanks to \eqref{eq:smalla} and \eqref{eq:z}, we can write
$$
\|\nabla\bigl(k(a)(\nabla\cA_\rho^{-1}\nabla a^\ell)\bigr)^\ell\|_{\dot B^s_{2,1}}\lesssim 
\|k(a)\|_{\dot B^{\frac d2}_{2,1}}\|\nabla\cA_\rho^{-1}\nabla a^\ell\|_{\dot B^{s+1}_{2,1}}\lesssim 
\|a\|_{\dot B^{\frac d2}_{2,1}}\| a\|_{\dot B^{s+1}_{2,1}}^\ell,$$
$$\begin{aligned}
\|\nabla\bigl(k(a)(\nabla\cA_\rho^{-1}\nabla \bigl(a^h +ak(a)))^\ell\bigr)\|_{\dot B^s_{2,1}}&\lesssim 
\|k(a)\|_{\dot B^{\frac d2}_{2,1}}\|\nabla\cA_\rho^{-1}\nabla (a^h +ak(a))\|^\ell_{\dot B^{s+1}_{2,1}}\\
&\lesssim 
2^{j_0}\|k(a)\|_{\dot B^{\frac d2}_{2,1}}\|\nabla\cA_\rho^{-1}\nabla  (a^h +ak(a))\|_{\dot B^{s}_{2,1}}\\&\lesssim 
\|a\|_{\dot B^{\frac d2}_{2,1}}\bigl(\| a\|_{\dot B^{s}_{2,1}}^h+
\|a\|_{\dot B^{\frac d2}_{2,1}}\| a\|_{\dot B^{s}_{2,1}}\bigr)\cdotp\end{aligned}$$
Hence, putting together with \eqref{eq:uw}, we get
\begin{equation}\label{eq:HF3v}
\|(\cA_1-\cA_\rho)w\|_{\dot B^s_{2,1}}\lesssim  \|a\|_{\dot B^{\frac d2}_{2,1}}\bigl(\|u\|^\ell_{\dot B^{s+2}_{2,1}}
+\|w\|^h_{\dot B^{s+2}_{2,1}}+\| a\|_{\dot B^{s+1}_{2,1}}^\ell+
\| a\|_{\dot B^{s}_{2,1}}^h+
\|a\|_{\dot B^{\frac d2}_{2,1}}\| a\|_{\dot B^{s}_{2,1}}\bigr)\cdotp
\end{equation}
Next, we have
$$\partial_t(\cA_\rho^{-1}\nabla Q(a))=\cA_\rho^{-1}\nabla (\d_t(Q(a)))+\partial_t(\cA_\rho^{-1})(\nabla Q(a)).$$
Computing $\d_t(Q'(a))$ from the mass equation, we see that
$$\cA_\rho^{-1}\nabla (\d_t(Q(a)))=\underbrace{-\cA_1^{-1}\nabla(Q'(a)\div((1+a)u))}_{I_1}
+\underbrace{(\cA_1^{-1}-\cA_\rho^{-1})\nabla(Q'(a)\div((1+a)u))}_{I_2}.$$
The term $I_1$  can be further decomposed by means of \eqref{eq:decompo}, and 
has already been treated  in the constant viscosity case. 
For $I_2,$ since 
$$\cA^{-1}_1-\cA_\rho^{-1}=\cA_1^{-1}(\cA_\rho-\cA_1)\cA_\rho^{-1},$$
we can write that 
$$I_2\simeq \cA_1^{-1}\nabla(k(a)\nabla\cA_\rho^{-1}\nabla(Q'(a)\div((1+a)u))).$$
Hence, using \eqref{eq:smalla}, \eqref{eq:z},  that $\nabla\cA_1^{-1}\nabla$ 
is homogeneous of degree $0$ and that
$s+1\leq d/2,$
\begin{equation}
\|I_2\|_{\dot B^s_{2,1}}^h
\lesssim 2^{-j_0}\|k(a)\nabla\cA_\rho^{-1}\nabla(Q'(a)\div((1+a)u))\bigr\|_{\dot B^s_{2,1}}\\
\lesssim \|a\|_{\dot B^{\frac d2}_{2,1}}\|u\|_{\dot B^{s+1}_{2,1}}.\end{equation}
Finally, a direct computation reveals that
\begin{equation*}
\partial_t(\cA_\rho^{-1})(\nabla Q(a))=\cA_\rho^{-1}\Bigl(\nabla\bigl(2\d_t\rho\,\mu'(\rho)D(\cA_\rho^{-1}\nabla(Q(a)))\bigr)
+\nabla\bigl(\d_t\rho\,\lambda'(\rho)\div (\cA_\rho^{-1}\nabla(Q(a)))\bigr)\Bigr)\cdotp\end{equation*}
Hence, 
$$
\partial_t(\cA_\rho^{-1})(\nabla Q(a))\simeq \cA_\rho^{-1}\nabla\Bigl((1+k(a))\div\bigl((1+a)u\bigr)\nabla\cA_\rho^{-1}\nabla(Q(a))\Bigr),$$
and we thus have
$$
\|\partial_t(\cA_\rho^{-1})(\nabla Q(a))\|_{\dot B^s_{2,1}}^h\lesssim 
2^{-j_0}\|\nabla\partial_t(\cA_\rho^{-1})(\nabla Q(a))\|_{\dot B^s_{2,1}}^h\lesssim
 \|a\|_{\dot B^{\frac d2}_{2,1}}\|u\|_{\dot B^{s+1}_{2,1}}.$$
 Compared to the constant viscosity case, all the new  terms in the right-hand side of \eqref{eq:wv}
 can be bounded by $\|a\|_{\dot B^{\frac d2}_{2,1}}\|u\|_{\dot B^{s+1}_{2,1}}.$
We conclude that Inequality \eqref{eq:wh} is still valid,  whence also  \eqref{eq:aw2}, and thus \eqref{eq:awnu} 
after rescaling. 
The inequalities corresponding to the limit case $s=-d/2$ can be adapted in the same way.

\subsection{Solving the limit equation}
Setting $q=1+b,$  Equation \eqref{eq:q} rewrites:
$$(\d_t+v\cdot\nabla)b+c^2(1+b)Q(b)=c^2(1+b)\div H_b\with 
H_b:=\bigl(\cA^{-1}_q-\cA_1^{-1}\bigr)(\nabla(Q(b)).$$
Compared to the constant coefficient case, one has to handle the nonzero right-hand side. 
Now, we have  the relation
$$\begin{aligned}
H_b&=\cA_q^{-1}(\cA_1-\cA_q)\cA_1^{-1}(\nabla(Q(b))\\
&=\cA_q^{-1} \!\Bigl(2\div\bigl((\mu(q)-\mu(1))D\bigl((-\Delta)^{-1}\nabla Q(b)\bigr)\bigr)\!+\!\cA_q^{-1}\nabla\bigl((\lambda(q)-\lambda(1)\bigr)\div\bigl((-\Delta)^{-1}\nabla Q(b)\bigr)\bigr)\Bigr)\cdotp
\end{aligned}$$
Hence $H_b$ is a combination of terms of the form $\cA_q^{-1}\nabla\bigl(k(b)\nabla(-\Delta)^{-1} \nabla Q(b) $
with $k(0)=0.$ 
As $\|b\|_{\dot B^{\frac d2}_{2,1}}$ is small, using \eqref{eq:z} yields
\begin{align}\label{eq:Hb}
\|(1+b)\div H_b\|_{\dot B^s_{2,1}}
&\lesssim  \|b\|_{\dot B^{\frac d2}_{2,1}}\|Q(b)\|_{\dot B^{s}_{2,1}}\nonumber\\
&\lesssim  \|b\|_{\dot B^{\frac d2}_{2,1}}\|b\|_{\dot B^{s}_{2,1}}.\end{align}
From this point, one can solve the limit equation as in the constant coefficient case whenever $b_0$ is small enough in $\dot B^{\frac d2}_{2,1}$
and get, for all $s\in(-d/2,d/2]$ and $t\in\R_+,$ 
$$\|b\|_{\wt L^\infty_t(\dot B^s_{2,1})} +\frac{c^2}2\int_0^t\|b\|_{\dot B^{s}_{2,1}}\leq \|b_0\|_{\dot B^{s}_{2,1}}.$$

\subsection{Passing to the limit}

The equation for $\da:=\ca-b$ now reads
 \begin{multline}\label{eq:dabis} \d_t\da
  + \check u\cdot\nabla\da +c^2(F(\check a)-F(b))=-(1+\check a)\div \check w +     (v-\check u)\cdot\nabla b
  \\+c^2\da\,\div H_b+c^2(1+\ca)\div(H_{\ca}-H_b).\end{multline}
Compared to \eqref{eq:da1}, only the last two terms are new. For the first one, we readily have by \eqref{eq:Hb}, 
$$\|\da\,\div H_b\|_{\dot B^{\frac d2-1}_{2,1}} \lesssim \|\da\|_{\dot B^{\frac d2-1}_{2,1}}\|b\|_{\dot B^{\frac d2}_{2,1}}^2.$$
For the last one, we use the following decomposition (where we omit the terms
pertaining to the shear viscosity  $\mu$ that can we treated exactly in the same way):
$$\begin{aligned}
H_{\ca}-H_b=(-\Delta)^{-1}\div\bigl(\cI_1+\cI_2+\cI_3)&\with\\
\cI_1&:=
(\lambda(\crho)-\lambda(q))\div(\cA_{\crho}^{-1}(\nabla Q(\ca))),\\
\cI_2&:=(\lambda(q)-\lambda(1))\div\bigl((\cA_{\crho}^{-1}-\cA_q^{-1})(\nabla Q(\ca))\bigr),\\
\andf\cI_3&:=(\lambda(q)-\lambda(1))\div\cA_q^{-1}\bigl(\nabla(Q(\ca)-Q(b))\bigr)\cdotp\end{aligned}$$
Since \eqref{eq:smalla} holds true, and  we want to bound  $(1+\ca)\div(H_{\ca}-H_b)$ in $L^1(\R_+;\dot B^{\frac d2-1}_{2,1}),$
it suffices to estimate $\cI_1,$ $\cI_2$ and $\cI_3$ in this space.
Now, owing to \eqref{eq:z}, we have
$$\begin{aligned}
\|\cI_1\|_{\dot B^{\frac d2-1}_{2,1}}&\lesssim \|(\lambda(\crho)-\lambda(q))\|_{\dot B^{\frac d2-1}_{2,1}}
\|\nabla(\cA_{\crho}^{-1}(\nabla Q(\ca)))\|_{\dot B^{\frac d2}_{2,1}}\\
&\lesssim
\|\da\|_{\dot B^{\frac d2-1}_{2,1}}\|\ca\|_{\dot B^{\frac d2}_{2,1}},\\[1ex]
\|\cI_3\|_{\dot B^{\frac d2-1}_{2,1}}&\lesssim \|(\lambda(q)-\lambda(1))\|_{\dot B^{\frac d2}_{2,1}}
\|\nabla\cA_q^{-1}(\nabla(Q(\ca)-Q(b)))\|_{\dot B^{\frac d2-1}_{2,1}}\\
&\lesssim \|b\|_{\dot B^{\frac d2}_{2,1}}\|Q(\ca)-Q(b)\|_{\dot B^{\frac d2-1}_{2,1}}\\&\lesssim
 \|b\|_{\dot B^{\frac d2}_{2,1}}\|\da\|_{\dot B^{\frac d2-1}_{2,1}}.\end{aligned}$$
To handle $\cI_2,$ we use the relation 
$$\cA_{\crho}^{-1}-\cA_q^{-1}=\cA_{\crho}^{-1}(\cA_q-\cA_{\crho})\cA_q^{-1}.$$  Hence, we have
$$
\|\cI_2\|_{\dot B^{\frac d2-1}_{2,1}}\lesssim \|b\|_{\dot B^{\frac d2}_{2,1}}
\|\nabla\cA_{\crho}^{-1}(\cA_q-\cA_{\crho})\cA_q^{-1}\nabla(Q(\ca))\|_{\dot B^{\frac d2-1}_{2,1}}.$$
Since 
$$
\cA_{\crho}-\cA_{q}=2\div\bigl((\mu(q)-\mu(\crho))D(\cdot)\bigr)+\div\bigl((\lambda(q)-\lambda(\crho))\nabla \cdot\bigr),
$$
it is easy to conclude, after taking advantage once more of \eqref{eq:z}, that
$$\|\cI_2\|_{\dot B^{\frac d2-1}_{2,1}}\lesssim \|b\|_{\dot B^{\frac d2}_{2,1}} \|\ca\|_{\dot B^{\frac d2}_{2,1}}
\|\da\|_{\dot B^{\frac d2-1}_{2,1}}.$$
To conclude, we have 
$$\|\da\,\div H_b+(1+\ca)\div(H_{\ca}-H_b)\|_{\dot B^{\frac d2-1}_{2,1}}
\lesssim \Bigl(\bigl(1+\|\ca\|_{\dot B^{\frac d2}_{2,1}}+\|b\|_{\dot B^{\frac d2}_{2,1}}\bigr)\|b\|_{\dot B^{\frac d2}_{2,1}}
+\|\ca\|_{\dot B^{\frac d2}_{2,1}}\Bigr)\|\da\|_{\dot B^{\frac d2-1}_{2,1}}.$$
Hence, as $\ca$ and $b$ are small in $\dot B^{\frac d2}_{2,1},$ one can absorb all the new terms 
created by variable viscosity  by the left-hand side of \eqref{eq:da3}, and conclude 
again that \eqref{eq:da5} holds true.

    \end{document}